\documentclass[leqno,11pt]{amsart}
\usepackage{amssymb,amsmath,latexsym,amsfonts,amsbsy,amsthm,mathtools,graphicx,color,subeqnarray,cases,tikz} 
\usepackage{amsfonts}
\usepackage{mathrsfs}
\usepackage{color}
\usepackage{amssymb,amsthm,amsmath,hyperref}
\usepackage{epsfig}
\usepackage{graphicx,float}
\usepackage{CJK}
\usepackage{setspace}
\usepackage{cancel}
\usepackage{verbatim}
\usepackage{ulem}
\usepackage{cite}

\usepackage{amssymb}
\newtheorem{theorem}{Theorem}[section]
\newtheorem{lemma}{Lemma}[section]
\newtheorem{remark}{Remark}[section]
\newtheorem{proposition}{Proposition}[section]

\newcommand{\beq}{\begin{eqnarray}}
\newcommand{\eeq}{\end{eqnarray}}
\newcommand{\beqno}{\begin{eqnarray*}}
\newcommand{\eeqno}{\end{eqnarray*}}
\newcommand{\be}{\begin{equation}}
\newcommand{\ee}{\end{equation}}
\newcommand{\beno}{\begin{equation*}}
\newcommand{\eeno}{\end{equation*}}

\def\e{{\rm e}}
\def\dd{{\rm d}}

\numberwithin{equation}{section}  
\allowdisplaybreaks[4]

\begin{document}
\title[Two-dimensional Poiseuille flow]{Nonlinear Stability of Two-Dimensional Poiseuille Flow via a Dilation Multiplier}

\author{Nangao Zhang }
\address[N. Zhang]{Corresponding Author, Research Center for Applied Mathematics and Interdisciplinary Sciences,  School of Mathematics and Statistics, Wuhan Textile University, Wuhan 430200, P. R. China.}
\email{mathngzhang@163.com}


\begin{abstract}

We investigate nonlinear stability of the Poiseuille flow \((y^2,0)\) for 2D incompressible Navier–Stokes equations on \(\mathbb T_x\times\mathbb R_y\).   We prove that if the initial vorticity satisfies \(\|\omega_{\mathrm{in}}\|_{H_x^1L_y^2}\ll\nu^{5/8}\) and the initial horizontal velocity perturbation satisfies \(\|P_0u^1_{\mathrm{in}}\|_{H_y^8}\ll 1\),  then the Poiseuille flow is globally nonlinearly stable. 
The key observation is that, after Fourier transformation, the Poiseuille flow and its curvature give rise to a kinetic operator and a potential whose commutators with the same bounded dilation multiplier are both nonnegative. 
We construct this multiplier as \(M=\arctan(A/4)\), where \(A\) generates dilations in the vertical Fourier variable, and prove that the sum of the two commutators controls the velocity. 
Then we use the modified energy to get inviscid damping and enhanced dissipation estimates. 
To close the nonlinear estimates, we incorporate the zero mode of the initial velocity into the background shear flow, thereby further improving the threshold.
\end{abstract}	

   
\maketitle
	
\tableofcontents

\section{Introduction}

We study the two-dimensional (2D) Navier--Stokes equations  in the domain $\Omega:=\mathbb{T}\times \mathbb{R}$:
\begin{equation}\label{1.1}
    \begin{cases}
\partial_t v + (v \cdot \nabla)v + \nabla P - \nu \Delta v = 0, \\[2mm]
\nabla \cdot v = 0,  \\[2mm]
v(0,x,y)=v_0(x,y), 
\end{cases}
\end{equation}
where \(\nu = \mathfrak{Re}^{-1} > 0\) is the viscosity coefficient and \(\mathfrak{Re}\) is the Reynolds number, \(v(t, x, y)= (v_1(t, x, y), v_2(t, x, y))\) is the velocity fields and \(P \) is the pressure.

It is easy to see
that the so-called Poiseuille flow:
\[
v_s(x,y)=(y^2,0),\ \ \ P_s=2\nu x,
\]
is a stationary solution of \eqref{1.1}.

Let's introduce the perturbation as
\[
u=v-v_s=(u_1,u_2),\ \ \ p=P-P_s,
\]
 then the perturbed equations around the Poiseuille
flow read as
\begin{equation}\label{1.2}
    \begin{cases}
\partial_t u - \nu \Delta u + (u \cdot \nabla)u + y^2\partial_x u + \begin{pmatrix} 2y u_2 \\ 0 \end{pmatrix} + \nabla p = 0, \\
\nabla \cdot u = 0, \\[2mm]
u(0,x,y) = u_{\mathrm{in}}(x, y).
\end{cases}
\end{equation}
Defining the vorticity of \(u\) as \(\omega := \nabla^\perp \cdot u\), where \(\nabla^\perp = (-\partial_y, \partial_x)\), it is possible to remove the pressure term and to rewrite the above system as
\begin{equation}\label{1.3}
    \begin{cases}
\partial_t \omega + y^2 \partial_x \omega - 2\partial_x \psi - \nu \Delta \omega = -u \cdot \nabla \omega, \\[2mm]
\omega = \Delta \psi,\ \ \  u = \nabla^\perp \psi,\\[2mm]
\omega(0,x,y)=\omega_{\mathrm{in}}(x,y).
\end{cases}
\end{equation}

Plane Poiseuille flow describes the steady motion of a viscous fluid driven by a constant pressure gradient between two stationary parallel plates. 
The pressure force balances viscous friction, producing a parabolic velocity profile.
  The quadratic profile \(U(y)=y^2\) studied here retains the vanishing shear and nonzero curvature at the centerline, while the absence of walls allows us to isolate their influence on the evolution of perturbations.
  
In this paper, we focus on quantifying a stability threshold in Sobolev spaces.  
Following \cite{BPM}, the problem can be formulated as follows: given a norm \(\|\cdot\|_X\) and find an exponent \(\gamma=\gamma(X)\) such that
\[
\|u_{{\rm in}}\|_X \lesssim \nu^{\gamma} \Longrightarrow \text{stability},
\]
\[
\ \ \ \ \ \ \ \ \ \ \ \ \   \|u_{{\rm in}}\|_X \gg \nu^{\gamma} \Longrightarrow \text{possible instability}.
\]
The exponent $\gamma$ is referred to as the {\it transition threshold}.

In recent years,  the stability of the Couette flow for the incompressible Navier–Stokes equations has attracted considerable interest among researchers, and substantial progress has been made.
More precisely, for the case of $(x,y)\in \mathbb{T}\times \mathbb{R}$,  it has been proved that
\begin{itemize}
    \item if $X$ is Gevrey class, then \(\gamma =0\) \cite{BMV},
    \item if $X$ is  Sobolev space, then \(\gamma \leq \frac{1}{2}\) \cite{BVW,Masmoudi2020},
    \item if $X$ is Sobolev space, then \(\gamma \leq \frac{1}{3}\) \cite{Masmoudi2022,WeiZhang2023},
    \item if \(X\) is Gevrey class \(\frac{1}{s}\), then \(\gamma \leq \frac{1-2s}{3(1-s)}\) \cite{LiZhao2025}. 
\end{itemize}
For the finite channel $(x,y)\in \mathbb{T}\times [-1,1]$, 
\begin{itemize}
    \item if $X$ is Sobolev space and the non-slip boundary condition is imposed, then \(\gamma\leq \frac{1}{2}\) \cite{CLWZ2020},
    \item if $X$ is Sobolev space and the Navier-slip boundary condition is imposed, then \(\gamma \leq \frac{1}{3}\) \cite{wei-zhang2026}.
\end{itemize}
For the infinite channel $(x,y)\in \mathbb{R}\times [-1,1]$, 
\begin{itemize}
    \item if $X$ is Sobolev space and the Navier-slip boundary condition is imposed, then \(\gamma\leq \frac{1}{2}+\) \cite{Arbon2025},
    \item if $X$ is Sobolev space and the Navier-slip boundary condition is imposed, then \(\gamma \leq \frac{1}{3}\) \cite{miao2025}.
\end{itemize}

For plane Poiseuille flow in \(\mathbb T\times[-1,1]\) with Navier-slip boundary conditions, Ding and Lin \cite{DingLin2022,DingLin2025} used resolvent estimates to obtain \(\gamma=2/3\) in the Sobolev
space. 
For the fixed-period channel with no-slip boundary conditions, Chen, Li, Shen, and Zhang \cite{ChenLi2025} established inviscid damping and enhanced dissipation for a class of symmetric flows, and also obtained \(\gamma=2/3\).
On \(\mathbb T\times[-1,1]\times\mathbb T\) with no-slip boundary condition, Chen, Ding, Lin, and Zhang \cite{Chen2025} proved nonlinear stability and obtained \(\gamma=7/4\) in the Sobolev
space. 

For 2D plane Poiseuille flow on \(\mathbb T\times\mathbb R\), Coti Zelati, Elgindi, and Widmayer \cite{Coti2020} used the hypocoercive energy to get nonlinear stability with \(\gamma=3/4+\). 
Del Zotto \cite{Del2023} removed the logarithmic loss from the linear decay time and obtained nonlinear enhanced dissipation under the condition
\[
\|\omega_{\rm in}\|_{L^2}
+\|P_0u_{\rm in}\|_{L^2}
\ll\,\frac{\nu^{2/3}}{1+|\log\nu|^{1/2}},
\]
where \(P_0\) denotes the streamwise average. 
More recently, Liang, Zhai, and Zhai \cite{Liang2026} obtained a log-free \(O(\nu^{2/3})\) sufficient scale in a viscosity-dependent weighted Sobolev space. 
For other results, see \cite{miao2025b,LiZhao2026}.

Based on the above observations, a natural question arises: {\bf can the stability threshold $2/3$ on the domain 
\(\mathbb T\times\mathbb R\) be further improved?} What is the key difficulty in achieving this goal? 
In this paper, we will give an affirmative answer to this question.

{\bf Notations.}
\begin{enumerate}
    \item[1.] For a function $f=f(t,x, y): \Omega \rightarrow \mathbb{R}$,  we define the Fourier transform as
\begin{align}\notag
\hat{f}_k(\eta) = \mathcal{F}(f)_k(\eta) = \frac{1}{2\pi} \iint_{\mathbb{T} \times \mathbb{R}} \mathrm{e}^{-i(kx+\eta y)} f(x, y) \,\mathrm{d}x\,\mathrm{d}y,\ \ \ \ k\in \mathbb{Z},\notag
\end{align}
and the inverse Fourier transform as
\[
\mathcal{F}^{-1}(\hat{f})(x, y) = \frac{1}{2\pi} \sum_{k \in \mathbb{Z}} \int_{\mathbb{R}} \mathrm{e}^{i(kx+\eta y)} \hat{f}_k(\eta) \,\mathrm{d}\eta.
\]
 \item[2.] For a function $f$, we define
 \[
P_0f(y)=\bar{f}(y) = \frac{1}{2\pi} \int_{\mathbb{T}} f(x,y) \,\mathrm{d}x, \quad P_{\neq} f =f_{\neq} =  f - \bar{f}, \quad f = \sum_{k \in \mathbb{Z}} f_k(y) \e^{ikx},
\]
where \(f_k(y)\) is the Fourier transform of $f$ in $x$.
\item[3.] When $1\leq p,q \leq \infty$ and $T_0\leq t\leq T_1$, we denote
\begin{align}\notag
& \|f\|_{L^q_{y}}=\|f\|_{L^q_{y}(\mathbb{R})},\ \ \ \ \ \ \ \|f\|_{L^q_{x,y}}=\|f\|_{L^q_{x,y}(\Omega)} ,\notag\\
&\|f\|_{L^pL^q}=\|f\|_{L_t^p(T_0,T_1)L^q_{x,y}(\Omega)},\ \ \ \ \|f\|_{H_{k,y}^{-1}}
:=
\|(-\partial_y^2+|k|^2)^{-\frac{1}{2}}f\|_{L_y^2}.\notag
\end{align}
The $H^q$-norm can be defined analogously.
\item[4.] We use the notation $a \lesssim b$  to indicate that  there exists a constant $C>0$ (independent of $\nu, t$) such that $a \leq C b$ (we define $a \gtrsim b$ analogously). 
In addition, we use the notation $a \approx b$ when there exists $C>0$ such that $C^{-1} b \leq a \leq C b$.
\item[5.] For given functions $f,g$, $1\leq q \leq \infty$, we denote
\begin{equation}\notag
 \|(f,g)\|_{L^q}=\|f\|_{L^q}+\|g\|_{L^q}, \ \ \ \  \langle f\rangle=\left(1+|f|^2\right)^{\frac{1}{2}}.
 \end{equation}
\end{enumerate}

Our main result is stated as follows.

\begin{theorem}\label{thm1}
   There exist a small constant \(c_*> 0\) ( independent of \(\nu\) and \(t\) ) such that if the initial vorticity \(\omega_{\mathrm{in}}\) and velocity \(P_0 u_{1,\mathrm{in}}\) satisfy
\begin{equation}\label{1.4}
\|P_0 u_{1,\mathrm{in}}\|_{H^6_y}\leq c_*,\ \ \ \ \varepsilon_* :=\|\omega_{\mathrm{in}}\|_{H^1_xL^2_{y}} \le c_* \nu^{\frac{5}{8}},
\end{equation}
then  there exists a small $\delta>0$ sufficiently small (depending only on universal constants) such that for all \(0 < \nu \le 1\) and \(0\leq t <\infty\), the perturbed equation \eqref{1.3} is global in time and satisfies the following stability estimates: 
\begin{equation}\label{1.5g}
    \sum_{k \in \mathbb{Z}} E_k \lesssim \varepsilon_*,
\end{equation}
where \(E_0=\|\overline{\omega}\|_{L^\infty L^2}\) and for \(k \neq 0\),
\begin{align}\label{1.5}
   E_k=& \|\e^{\delta\nu^{\frac{1}{2}}t}\omega_k\|_{L^\infty L^2} + \nu^{\frac{1}{2}}\|\e^{\delta\nu^{\frac{1}{2}}t}(k,\partial_y)\omega_k\|_{L^2 L^2} + (\nu|k|)^{\frac{1}{4}}\|\e^{\delta\nu^{\frac{1}{2}}t}\omega_k\|_{L^2 L^2}  \notag\\
& + |k|^{\frac{1}{2}}\|\e^{\delta\nu^{\frac{1}{2}}t}u_k\|_{L^2 L^2}.
\end{align}

\end{theorem}
\begin{remark}
In fact, from \eqref{1.5g}--\eqref{1.5}, we can derive that
\[
\|\omega_{\ne}\|_{L^2}
\lesssim \varepsilon_*\e^{-\delta\nu^{\frac{1}{2}}t},
\qquad
\bigl\||D_x|^{1/2}u_{\ne}\bigr\|_{L^2L^2}
\lesssim \varepsilon_*\e^{-\delta\nu^{\frac{1}{2}}t}.
\]
 These two estimates correspond to enhanced dissipation and the time-integrated form of inviscid damping, respectively.
\end{remark}
\begin{remark}
The exponent \(1\) in the hypothesis \(H_x^1L_y^2\) is not essential.  For every fixed \(s>1/2\), Cauchy–Schwarz gives
\[
\sum_{k\in\mathbb Z}\|\omega_{\mathrm{in},k}\|_{L_y^2}
\leq
\left(\sum_{k\in\mathbb Z}(1+k^2)^{-s}\right)^{1/2}
\|\omega_{\mathrm{in}}\|_{H_x^sL_y^2}.
\]
Hence the same proof applies under the assumption \(\|\omega_{\mathrm{in}}\|_{H_x^sL_y^2}\leq c_s\nu^{5/8}\), with the condition on \(P_0u_{1,\mathrm{in}}\) unchanged. 
\end{remark}
\begin{remark}
    The condition on \(V_{\rm in}:=P_0u_{1,\rm in}\) ensures that its heat evolution \(V^L=\e^{\nu t\partial_y^2}V_{\rm in}\) can be incorporated into the background shear flow.
\end{remark}
\begin{remark}
    We believe that this results can also be extended to the 2D Navier Stokes Boussinesq equations, where the threshold for vorticity remains unchanged, while the threshold for temperature is expected to be $9/8$.
\end{remark}

\section{Outline of the proof}\label{s2}

\subsection{Refined decomposition}

In this subsection, we refine the decomposition of the perturbed equation \eqref{1.3} according to its structure. These insights are essential for deriving the improved stability threshold subsequently.

First, we denote \(V_{\mathrm{in}}(y) = P_0u_{1,\mathrm{in}}(y)\), and define the following function:
\begin{equation}\label{2.1}
V^L(t,y) = \e^{\nu t \partial_y^2} V_{\mathrm{in}}(y). 
\end{equation}
It is obvious that $V^L(t,y)$ satisfies equation 
\begin{equation}\label{2.2}
    \partial_tV^L-\nu\partial_y^2V^L=0,\ \ \ \ V^L(0,y)=V_{\mathrm{in}}(y).
\end{equation}
 This is actually the linearized equation satisfied by for $\bar{u}_1(t,y)$. 
 Following this, we introduce 
 \begin{equation}\label{2.3}
   U(t,y):=\bar{u}_1-V^L , \ \ \  \ \ V(t,y) := y^2 + V^L(t,y),
 \end{equation}
 then it holds that
\begin{align}
    \omega&=\bar{\omega}+\omega_{\neq}=-\partial_y\bar{u}_1+\omega_{\neq}=-\partial_yV^L-\partial_yU+\omega_{\neq},\notag\\
    u&=\bar{u}+u_{\neq} = (V^L + U, 0) +u_{\neq},\notag\\
    u\cdot \nabla \omega&=[(V^L + U, 0) +u_{\neq}]\cdot\nabla[-\partial_yV^L-\partial_yU+\omega_{\neq}]\notag\\
    &=(V^L + U)\partial_x\omega_{\neq}-u_{2,\neq}\partial_y^2V^L+u_{\neq}\cdot\nabla(-\partial_yU+\omega_{\neq}).\notag
\end{align}
Therefore, we can rewrite \eqref{1.3} as
\begin{equation}\label{2.4}
    \begin{cases}
\partial_t\omega-\nu \Delta\omega+ V\partial_x\omega- \partial_y^2V\partial_x \Delta^{-1}\omega=-\nabla \cdot [u_{\neq}(\omega_{\neq}-\partial_yU) + (U\omega_{\neq}, 0)],\\[2mm]
U_t - \nu \partial_y^2U = \overline{u_{2,\neq} \omega_{\neq}},\\[2mm]
 \omega|_{t=0} =\omega_{\mathrm{in}}(x,y), \ \ U|_{t=0} =0.     
    \end{cases}
\end{equation}

\subsection{Key ideas}

In this section, we provide a detailed explanation of the ideas involved in the proofs throughout this paper. 
For the convenience of handling the nonlinear system later, we first consider a linear system after performing the Fourier transform only with respect to the variable $x$:
\begin{equation}\label{2.5}
\begin{cases}
\partial_tf_{k} + \nu \Delta_k f_{k} + i k(y^2 + 2\Delta_k^{-1}) f_{k} = \mathrm{div}_k F+\nu G_1-ikG_2+G_3, \quad\\[2mm] 
f_{k}(0,y) = f_{k, \mathrm{in}},\ \ \ k\in \mathbb Z~ \mbox{and}~k\neq 0,~y\in \mathbb{R},
\end{cases}
\end{equation}
where \(\Delta_k = -\partial_y^2 + k^2\), \(\mathrm{div}_k F = i k F_1 + \partial_y F_2\).
Then taking the Fourier transform in \(y\), we have from \eqref{2.5} that
\begin{equation}\label{2.6}
\partial_t\hat{f}_k + \nu (|k|^2 + \eta^2) \hat{f}_k + i k\left(-\partial_\eta^2+\frac{2}{|k|^2 + \eta^2}\right)  \hat{f}_k = i(k \hat{F_1} + \eta \hat{F_2})+\nu \hat G_1-ik\hat G_2+\hat G_3. 
\end{equation}

At this point, the nonlocal term \(\Delta_k^{-1}\) in physical space becomes an explicit multiplication potential in the $\eta$-space, while $y^2$ becomes second-order derivative. 
From this, we observe that the equation contains a Schrödinger operator with a special potential.
To handle the nonlinear terms later, we need to exploit some favorable mechanisms of the linearized equation, such as enhanced dissipation and inviscid damping.
Nevertheless, previous results indicate that deriving inviscid damping is rather challenging, since the Poiseuille flow is a non-monotonic flow.
Consequently, the threshold could not be further improved for a considerably long time.
Through careful analysis, we find that it is necessary to seek an additional energy correction such that the Schrödinger operator \(i \left(-\partial_\eta^2+\frac{2}{|k|^2 + \eta^2}\right)\), which disappears in the standard $L^2$ energy estimate, yields a positive term that controls \(\|(k^2+\eta^2)^{-1/2}\hat{f}_k\|_{L^2}\).
 This is the key to obtaining inviscid damping.
To overcome this obstacle, we select a dilation generator \(A\) (see \eqref{3.6b} in Lemma \ref{l1}), whose commutators with both the kinetic and potential of the Schrödinger operator are positive. 
That is to say, we can prove 
\[
i\left[-\partial_\eta^2+\frac{2}{|k|^2 + \eta^2},A\right]
=
-2\partial_\eta^2+\frac{4\eta^2}{(|k|^2+\eta^2)^2}
\ge0.
\]
However, since \(A\) is an unbounded operator, 
\[
|\langle Ah,h\rangle|\lesssim\|h\|_{L^2}^2
\]
 cannot be guaranteed. 
 If \(A\) is directly introduced into the modified energy, it is impossible to obtain the energy functional that is equivalent to the standard $L^2$ estimate. Therefore, we must next bound the dilation generator while preserving the positive commutator structure.
{\bf Through analysis, we find that the multiplier \(M=\arctan(A/4)\) can be chosen because it simultaneously satisfies boundedness, a favorable commutator structure, and compatibility with the weighted spaces}, see \eqref{eq:dilation-F-M-m}--\eqref{eq:dilation-M-basic}, \eqref{2.26d}--\eqref{2.27d}, \eqref{2.30d}--\eqref{2.31d}, \eqref{b} and \eqref{3.77c}.
With these preparations completed, we can obtain the estimates for the linearized equation \eqref{2.5}, which include those for inviscid damping and enhanced dissipation, see Proposition \ref{p3.2} for details.

Following the treatment of the linearized equation, we proceed to analyze nonlinear stability and the  transition threshold problem.
We find that for the nonlinear term \(u\cdot \nabla \omega\), \(\bar{u}_1\partial_x\omega\) is the worst term. 
In fact, if we treat this term within the nonlinear part, we can close our estimates under \(\|\omega_{\mathrm{in}}\|_{H^1_xL^2_{y}} \le c_* \nu^{\frac{2}{3}},\) this paper gives a
new proof different from the method introduced in \cite{Del2023}.
However, in that case, we would only remove the logarithmic loss in the threshold obtained in \cite{Del2023}, which is not a significant improvement.
To further improve the transition threshold, we have to introduce new ideas. 
{\bf Our second key idea is to treat the term given by the product of \(\partial_x\omega\) and the solution $V^L(t,y)$ to the linearized part of the equation for \(\bar{u}_1\) as a part of the linearized equation for \(\omega\), thereby introducing a new shear flow $V(t,y)$}. See \eqref{2.1}--\eqref{2.4} for more details.
It is worth noting that, since the scale of $U$ is larger here, it ultimately does not influence the final threshold.
However, we must pay a price for doing this: the structure of the decomposed linearized equation no longer matches those we have investigated before, and we are still unable to obtain the desired results. 
For this purpose, {\bf our third idea is to introduce a map \(\mathcal T_{k,t}\) that transforms the decomposed linear equation into the form of the equations we have studied, and subsequently estimate the error terms it generates}, see \eqref{4.6c}--\eqref{4.9cc}. 
Then, by using Proposition \ref{p3.2}, we finally obtain the estimates for the linearized equation of \eqref{2.4}.
Subsequently, by exploiting some favorable properties of the linear equations, we ultimately obtain the result of nonlinear stability and improve the transition threshold to $\frac{5}{8}$.

Finally, we further discuss the optimality of the threshold. 
Since the Poiseuille flow is non-monotonic and its enhanced dissipation rate is weaker compared to that of monotonic flows (such as Couette flow), it is natural to conjecture that the optimal threshold might be $\frac{1}{2}$. However, a recent paper \cite{ChenJia}, considering  the asymptotic stability of the Kolmogorov flow on a non-square torus, suggests that transition threshold $\frac{1}{3}$ might also be achievable for non-monotonic flows. 
In their work, they employed not only enhanced dissipation and inviscid damping, but also another favorable physical mechanism, namely vorticity depletion. 
Therefore, we have reason to believe that if one further derives the vorticity depletion estimate for the linearized equation based on this paper, and makes some more delicate estimate of the nonlinear system by using the good derivative  $\partial_y+2yt\partial_x$, improvement of the threshold is highly likely. This will be considered in our forthcoming work.

\section{Linear inviscid damping and enhanced dissipation for Poiseuille flow}\label{s3}

In this section, we devote ourselves to deriving the estimates for the linearized equation \eqref{2.5}.
More precisely, we first introduce a bounded multiplier $M$, and subsequently obtain the expected estimates through the construction of an appropriate energy functional.

\subsection{Construction of the unitary group and bounded multiplier}

We first introduce the unitary dilation group on \(L^2(\mathbb R)\) and identify its self-adjoint generator $A$.

\begin{lemma}\label{l1}
 On the complex Hilbert space
$L^2(\mathbb R_\eta)$, for any $h=h(\eta)\in L^2(\mathbb R)$, we define the operator
\begin{equation}\label{2.9b}
 (\mathcal{U}_s h)(\eta)=\e^{s/2}h(\e^s\eta),\qquad s\in\mathbb R ,
\end{equation}
then $(\mathcal{U}_s)_{s\in\mathbb R}$ is a strongly continuous unitary group.
Denote its self-adjoint generator by $A$, with the convention $\mathcal{U}_s=\e^{isA}$. We also define the map
\begin{align}\label{3.7b}
 J:L^2(\mathbb R)&\longrightarrow
       L^2(\mathbb R)\oplus L^2(\mathbb R),\qquad
       Jh=(J_+h,J_-h),\notag\\
 (J_\pm h)(z)&=\e^{z/2}h(\pm \e^z),
 \end{align}
then $J$ is unitary, and
\begin{align}\label{3.8b}
 D(A)&=\{h \in L^2(\mathbb{R}) : Jh \in H^1(\mathbb{R}) \oplus H^1(\mathbb{R})\},\notag\\
 JAJ^{-1}&=(-i\partial_z)\oplus(-i\partial_z).
 \end{align}
Here $\oplus$ denotes the direct sum, and its elements are pairs of functions. For every $h\in D(A)$,
\begin{equation}\label{3.6b}
 Ah=-i\left(\eta\,\partial_\eta h+\frac12 h\right).
\end{equation}
In addition, if $h\in L^2(\mathbb R_\eta)$ and $\eta\,\partial_\eta h\in L^2(\mathbb R_\eta)$, we can get $ h\in D(A)$.

Finally, let 
\begin{equation}\label{D}
  D(\eta):=\{h\in L^2(\mathbb R_\eta):\eta h\in L^2(\mathbb R_\eta)\}.  
\end{equation}
Then for every $s\in\mathbb R$,
\begin{equation}\label{3.6bb}
 \eta \mathcal{U}_s h=\e^{-s}\mathcal{U}_s(\eta h)\quad(h\in D(\eta)),\qquad
 -i\partial_\eta(\mathcal{U}_s h)=\e^s\mathcal{U}_s(-i\partial_{\eta}h)\quad(h\in H^1(\mathbb R_\eta)).
\end{equation}

\end{lemma}

\begin{proof}
 For any \(h \in L^2(\mathbb{R})\), by making the change of variables \(r = {\rm e}^s\eta\), we obtain
\[
\|\mathcal{U}_s h\|_{L^2_\eta}^2 = \int_{\mathbb{R}} {\rm e}^s |h({\rm e}^s\eta)|^2 \dd \eta = \int_{\mathbb{R}} |h(r)|^2 \dd r = \|h\|_{L^2_\eta}^2.
\]
Thus \(\mathcal{U}_s\) is an isometry. 
Naturally, we can also define \(\mathcal{U}_{-s}\). By computing the composite operator $\mathcal{U}_{-s}\mathcal{U}_s$, we get
\[
(\mathcal{U}_{-s}\mathcal{U}_s h)(\eta) = \mathcal{U}_{-s}(\e^{s/2}h(\e^s\eta)) = \e^{-s/2}[\e^{s/2}h(\e^{-s} \cdot \e^s\eta)] = h(\eta).
\]
Similarly, we can get \(\mathcal{U}_s \mathcal{U}_{-s} = I\), where $I$ is an identity operator.
Thus, \(\mathcal{U}_{-s}\) is the inverse of \(\mathcal{U}_s\), \(\mathcal{U}_s\) is a unitary operator.

In addition, for any \(s, t \in \mathbb{R}\),
\[
(\mathcal{U}_s(\mathcal{U}_t h))(\eta) = \e^{s/2}(\mathcal{U}_t h)(\e^s\eta) = \e^{s/2}\e^{t/2}h(\e^t(\e^s\eta)) = \e^{(s+t)/2}h(\e^{(s+t)}\eta)
\]
which is exactly equal to \((\mathcal{U}_{s+t}h)(\eta)\), namely,
\[
\mathcal{U}_s \mathcal{U}_t = \mathcal{U}_{s+t}.
\]
Setting \(s = 0\) and substituting it into the definition \eqref{2.9b}:
\[
(\mathcal{U}_0 h)(\eta) = \e^0 h(\e^0\eta) = h(\eta),
\]
i.e., \(\mathcal{U}_0 = I\). 
 Therefore, \((\mathcal{U}_s)_{s \in \mathbb{R}}\) is a unitary group.

To verify strong continuity, we first take \(h \in C_c^\infty(\mathbb{R})\). As \(s \to 0\),
\[
{\rm e}^{s/2} h({\rm e}^s\eta) \to h(\eta)
\]
pointwise. For \(s\) restricted to a neighborhood of zero, the support of the function \(\mathcal{U}_s h\) is contained in a fixed compact set, and it is uniformly bounded on this set. Then by Lebesgue's dominated convergence theorem,
\[
\|\mathcal{U}_s h - h\|_{L^2} \to 0.
\]
Since \(C_c^\infty(\mathbb{R})\) is dense in \(L^2(\mathbb{R})\) and \(\|\mathcal{U}_s\|_{L^2 \to L^2} = 1\), the above convergence extends to all \(h \in L^2(\mathbb{R})\). Hence the group is strongly continuous.

 For \(h \in L^2(\mathbb{R}_\eta)\), 
By using \eqref{3.7b} and introducing the changes of variables \(\eta = {\rm e}^z\) or \(\eta = -{\rm e}^z\), respectively, we have
\[
\int_{\mathbb{R}} |J_+ h(z)|^2 \dd z = \int_0^\infty |h(\eta)|^2 \dd \eta,
\]
and
\[
\int_{\mathbb{R}} |J_- h(z)|^2 \dd z = \int_{-\infty}^0 |h(\eta)|^2 \dd \eta.
\]
Therefore,
\begin{equation}\label{jh}
   \|J h\|_{L^2}^2=\|J_+ h\|_{L^2}^2 + \|J_- h\|_{L^2}^2 = \|h\|_{L^2}^2. 
\end{equation}
This means that $J$ is an isometry.
Conversely, for any \(g_\pm(z) \in L^2(\mathbb{R}_z)\), we can define separately on the positive and negative half-axes
\[
h(\eta) = \begin{cases} \eta^{-1/2} g_+(\log \eta), & \eta > 0, \\[2mm] |\eta|^{-1/2} g_-(\log |\eta|), & \eta < 0. \end{cases}
\]
For the positive half-axis, note that
$$
\int_0^\infty |h(\eta)|^2 \dd\eta = \int_0^\infty \eta^{-1} |g_+(\log \eta)|^2 \dd\eta.
$$
Let
$$
z = \log \eta, \quad \eta = \e^z, \quad \dd\eta = \e^z \dd z = \eta \dd z,
$$
then
$$
\eta^{-1} \dd\eta = \dd z,
$$
and 
$$
\int_0^\infty |h(\eta)|^2 \dd\eta = \int_{\mathbb{R}} |g_+(z)|^2 \dd z.
$$
Similarly, for the negative half-axis,
$$
\int_{-\infty}^0 |h(\eta)|^2 \dd\eta = \int_{\mathbb{R}} |g_-(z)|^2 \dd z.
$$
Therefore, we arrive at
$$
\|h\|_{L^2(\mathbb{R}_\eta)}^2 = \int_{\mathbb{R}} |g_+(z)|^2 \dd z + \int_{\mathbb{R}} |g_-(z)|^2 \dd z = \|(g_+, g_-)\|_{L^2(\mathbb{R}_z) }^2,
$$
and then \(h \in L^2(\mathbb{R}_\eta)\).
Therefore, for every \((g_+,g_-)\in L^2(\mathbb R_z)\oplus L^2(\mathbb R_z)\), the inverse construction above produces an \(h\in L^2(\mathbb R_\eta)\) satisfying
\begin{equation}\label{3.10aaa}
    (Jh)(z)=((J_+h)(z),(J_-h)(z)) = (g_+(z), g_-(z)).
\end{equation}
Thus, \(J\) is surjective and hence a unitary mapping from \(L^2(\mathbb{R}_\eta)\) to \(L^2(\mathbb{R}_z) \oplus L^2(\mathbb{R}_z)\).
In addition, for every \(h\in L^2(\mathbb R_\eta)\) and \(s \in \mathbb{R}\), a direct computation yields
\[
\begin{aligned}
(J_\pm \mathcal{U}_s h)(z) &= {\rm e}^{z/2} {\rm e}^{s/2} h(\pm {\rm e}^s {\rm e}^z) \\
&= {\rm e}^{(z+s)/2} h(\pm {\rm e}^{z+s}) \\
&= (J_\pm h)(z+s)=g_{\pm}(z+s):=(T_s g_{\pm})(z).
\end{aligned}
\]
Thus \(J \mathcal{U}_s J^{-1}\) consists of two identical translation groups, and
\begin{equation}\label{2.13b}
J(\mathcal{U}_s h) = (T_s g_+, T_s g_-).
\end{equation}
We now transfer this characterization to \(\mathcal{U}_s\).
Since \(\mathcal{U}_s=\e^{isA}\), the domain of its self-adjoint generator is characterized by
\[
h\in D(A)
\quad\Longleftrightarrow\quad
\lim_{s\to0}\frac{\mathcal{U}_sh-h}{s}
\text{ exists in }L^2(\mathbb R_\eta),
\]
and the limit is \(iAh\). By \eqref{3.10aaa}--\eqref{2.13b}, we infer
\begin{equation}\label{2.14b}
J\!\left(\frac{\mathcal{U}_sh-h}{s}\right)
=
\left(
\frac{T_sg_+-g_+}{s},
\frac{T_sg_--g_-}{s}
\right).
\end{equation}
Because \(J\) is unitary, the quotient on the left before applying \(J\) converges in \(L^2(\mathbb R_\eta)\) if and only if both quotients on the right converge in \(L^2(\mathbb R_z)\). The characterization of the scalar translation generator therefore gives
\[
h\in D(A)
\quad\Longleftrightarrow\quad
g_+\in H^1(\mathbb R_z)
\ \text{and}\
g_-\in H^1(\mathbb R_z).
\]
Finally, for \(h\in D(A)\), taking \(s\to0\) in the difference-quotient identity \eqref{2.14b} yields
\[
iJ(Ah)=(g_+',g_-').
\]
Equivalently,
\begin{equation}\label{abc}
J(Ah)=(-ig_+',-ig_-').
\end{equation}
Since \(Jh=(g_+,g_-)\) and \(h=J^{-1}(g_+,g_-)\), we conclude that
\[
(JAJ^{-1})(g_+,g_-)=J(Ah)=(-i\partial_zg_+,-i\partial_zg_-),
\]
with domain
\[
D(JAJ^{-1})
=H^1(\mathbb R_z)\oplus H^1(\mathbb R_z).
\]
Thus, we get \eqref{3.8b}.

 On both half-axes,  from \eqref{3.7b}, we have 
\[
\begin{aligned}
\partial_z (J_+ h)(z) &= \frac{1}{2} {\rm e}^{z/2} h({\rm e}^z) + {\rm e}^{z/2} {\rm e}^z h'({\rm e}^z) \\
&= J_+ \left( \frac{1}{2} h + \eta h' \right)(z).
\end{aligned}
\]
For the negative half-axis, since \(\partial_z h(-{\rm e}^z) = -{\rm e}^z h'(-{\rm e}^z)\) and \(\eta = -{\rm e}^z\), we still have
\[
\partial_z (J_- h) = J_- \left( \frac{1}{2} h + \eta h' \right).
\]
Therefore, it holds that

\[
\partial_z(J_\pm h)
   =J_\pm\!\left(\eta\,\partial_\eta h+\tfrac12h\right),
\]
which together with \eqref{abc} gives
\[
\begin{aligned}
J(Ah)
&=\bigl(-i\partial_z(J_+h),\,-i\partial_z(J_-h)\bigr)\\
&=\left(
J_+\!\left[-i\left(\eta\partial_\eta h+\tfrac12h\right)\right],
J_-\!\left[-i\left(\eta\partial_\eta h+\tfrac12h\right)\right]
\right)\\
&=J\!\left[-i\left(\eta\partial_\eta h+\tfrac12h\right)\right].
\end{aligned}
\]
Thus \eqref{3.6b} holds on the whole line.

Finally, for $h\in D(\eta)$, a direct computation from \eqref{2.9b} gives
\[
\begin{aligned}
\eta (\mathcal{U}_s h)(\eta) &= \eta {\rm e}^{s/2} h({\rm e}^s\eta) = {\rm e}^{-s} {\rm e}^{s/2} ({\rm e}^s\eta) h({\rm e}^s\eta) \\
&= {\rm e}^{-s} (\mathcal{U}_s(\eta h))(\eta).
\end{aligned}
\]
If \(h \in H^1(\mathbb{R})\), 
\[
\partial_\eta (\mathcal{U}_s h)(\eta) = {\rm e}^{3s/2} \partial_\eta h({\rm e}^s\eta),
\]
and then
\[
{\rm e}^s \mathcal{U}_s(-i\partial_\eta h)={\rm e}^{3s/2}(-i\partial_\eta h)({\rm e}^s\eta)=-i{\rm e}^{3s/2}(\partial_{\eta}h)({\rm e}^s\eta)=-i\partial_{\eta}(\mathcal{U}_s h).
\]
Therefore, we get \eqref{3.6bb}, and hence completes the proof of this lemma. 

\end{proof}

On the basis of Lemma \ref{l1}, we next define 
\begin{equation}\label{eq:dilation-F-M-m}
F(a)=\arctan(\frac{a}{4}),
\qquad
M=F(A),
\qquad
m(a)=\frac12\log\frac{a^2+25}{a^2+9}.
\end{equation}
Since $F(a)$ is real-valued and $|F(a)|\le\pi/2$, by Lemma \ref{l1} and Plancherel's formula,  one can verify that
\begin{equation}\label{eq:dilation-M-basic}
M=M^*,
\qquad
\|M\|_{L^2\to L^2}\le\frac{\pi}{2}.
\end{equation}
This means that 
$M$ is a bounded self-adjoint operator, even though
the generator $A$ is unbounded. 

The following lemma provides an integral representation suitable
for the subsequent commutator calculations.

\begin{lemma}
\label{l2}
For   $\epsilon>0$, $a\in\mathbb R$ and $h\in L^2(\mathbb R)$, define
\begin{equation}\label{eq:dilation-M-epsilon}
M_\epsilon h
=
\int_\epsilon^\infty
\frac{\e^{-4s}}{2is}
\bigl(\mathcal U_sh-\mathcal U_{-s}h\bigr)\dd s, \qquad
F_\epsilon(a)
=
\int_\epsilon^\infty
\e^{-4s}\frac{\sin(as)}s\dd s,
\end{equation}
where $\mathcal U_s$ is defined in \eqref{2.9b}. Then $M_\epsilon$ is a bounded self-adjoint operator and
\begin{equation}\label{eq:dilation-M-epsilon-uniform}
|F_\epsilon(a)|\lesssim 1,\ \ \  \ \  \|M_\epsilon\|_{L^2\to L^2}\lesssim 1.
\end{equation}
Moreover, it holds
\begin{equation}\label{eq:dilation-F-epsilon}
M_\epsilon=F_\epsilon(A),
\end{equation}
and
\begin{equation}\label{eq:dilation-M-strong}
\lim_{\epsilon\to0}
\|(M_\epsilon-M)h\|_{L^2}=0.
\end{equation}
The function $m$ satisfies
\begin{equation}\label{eq:dilation-m-bounds}
\frac4{a^2+16}
\le m(a)
\le\frac{16}{a^2+16},
\qquad a\in\mathbb R.
\end{equation}
Here $M$ and $m$ are defined in \eqref{eq:dilation-F-M-m}.
\end{lemma}

\begin{proof}
Let \(\mathcal F_z\) denote the Fourier transform in \(z\), with the convention
\[
(\mathcal F_zg)(a)
=\frac1{\sqrt{2\pi}}\int_{\mathbb R}\e^{-iaz}g(z)\dd z.
\]
Recall the definition of $J$ in \eqref{3.7b}, we now define the mapping $\mathcal V:
L^2(\mathbb R_\eta)
\longrightarrow
L^2(\mathbb R_a)\oplus L^2(\mathbb R_a)$ as
\begin{equation}\notag
 \mathcal V:=(\mathcal F_z\oplus\mathcal F_z)J.   
\end{equation}
That is, for any $h\in L^2(\mathbb R_\eta)$, 
\begin{equation}\label{V}
\mathcal Vh=(\mathcal F_z(J_+ h),\mathcal F_z(J_- h)).
\end{equation}
Since \(J\) and the direct sum \(\mathcal F_z\oplus\mathcal F_z\) are unitary, their composition \(\mathcal V=(\mathcal F_z\oplus\mathcal F_z)J\) is a unitary map from \(L^2(\mathbb R_\eta)\) onto \(L^2(\mathbb R_a)\oplus L^2(\mathbb R_a)\). In particular, we can prove
\[
\|\mathcal F_z(J_+ h)\|_{L^2_a}^2+\|\mathcal F_z(J_- h)\|_{L^2_a}^2
=\|h\|_{L^2_\eta}^2.
\]
We next compute the action of the dilation group under \(\mathcal V\). The identity \eqref{2.13b}
shows that dilation becomes translation after applying \(J\). Under the Fourier transform, translation satisfies
\[
\mathcal F_z[g(\,\cdot+s)](a)
=\e^{isa}\mathcal F_zg(a).
\]
 Applying this to both components of \(Jh\), we obtain
\begin{equation}\label{V1}
\mathcal V(\mathcal U_sh)
=(\e^{isa}\mathcal F_z(J_+ h),\e^{isa}\mathcal F_z(J_- h)),
\end{equation}
and, replacing \(s\) by \(-s\),
\begin{equation}\label{V2}
\mathcal V\mathcal U_{-s}h
=(\e^{-isa}\mathcal F_z(J_+ h)(a),\e^{-isa}\mathcal F_z(J_- h)(a)).
\end{equation}
Thus, by applying \eqref{eq:dilation-M-epsilon}, we infer
\[
\begin{aligned}
\mathcal VM_\epsilon h
&=\int_\epsilon^\infty
\frac{\e^{-4s}}{2is}
\bigl(\mathcal V\mathcal U_sh-\mathcal V\mathcal U_{-s}h\bigr)\dd s\\
&=\int_\epsilon^\infty
\frac{\e^{-4s}}{2is}
\left(
(\e^{isa}-\e^{-isa})\mathcal F_z(J_+ h),\,
(\e^{isa}-\e^{-isa})\mathcal F_z(J_- h)
\right)\dd s.
\end{aligned}
\]
Using 
\[
\frac{\e^{isa}-\e^{-isa}}{2i}=\sin(as),
\]
we therefore obtain
\[
\mathcal VM_\epsilon h
=\bigl(F_\epsilon(a)\mathcal F_z(J_+ h)(a),F_\epsilon(a)\mathcal F_z(J_- h)(a)\bigr).
\]
By \eqref{abc} and the definition of \(\mathcal V\), for \( h\in D(A)\), we have
\begin{equation}\label{ah}
\mathcal V(Ah)
=(\mathcal F_z\oplus\mathcal F_z)\bigl(J(Ah)\bigr)
=(\mathcal F_z(-ig_+'),\mathcal F_z(-ig_-'))
=(a\mathcal F_z(J_+ h),a\mathcal F_z(J_- h)).
\end{equation}
Consequently, the bounded spectral multiplier \(F_\epsilon(A)\) is defined by
\begin{equation}\label{F}
   F_\epsilon(A)h
=\mathcal V^{-1}(F_\epsilon \mathcal F_z(J_+ h),F_\epsilon \mathcal F_z(J_- h)). 
\end{equation}
Comparing this expression with the formula for \(\mathcal VM_\epsilon h\), we conclude that
\[
M_\epsilon=F_\epsilon(A).
\]
Since \(F_\epsilon\) is real-valued, multiplication by \(F_\epsilon\) is self-adjoint on each component. Unitary conjugation then gives \(M_\epsilon=M_\epsilon^*\). Thus, we obtain \eqref{eq:dilation-F-epsilon}.

We now establish a bound for \(F_\epsilon\) that is uniform in both \(\epsilon\) and \(a\). Clearly, \(F_\epsilon(0)=0\). For \(a\ne0\), introduce
\[
r=|a|s,\qquad b=\frac4{|a|},\qquad \rho=|a|\epsilon.
\]
Then
\[
\e^{-4s}=\e^{-br},
\qquad
\frac{\dd s}{s}=\frac{\dd r}{r},
\qquad
\sin(as)=\operatorname{sgn}(a)\sin r,
\]
and hence
\[
F_\epsilon(a)
=\operatorname{sgn}(a)
\int_\rho^\infty \e^{-br}\frac{\sin r}{r}\dd r.
\]
This substitution puts the oscillation into the fixed form \(\sin r\).

If \(\rho<1\), the part of the integral over \([\rho,1]\) is controlled by \(|\sin r|\le r\):
\[
\left|\int_\rho^1\e^{-br}\frac{\sin r}{r}\dd r\right|
\le\int_\rho^1 1\dd r
\le1.
\]
If \(\rho\geq 1\), set
\[
w(r)=\frac{\e^{-br}}r,\quad r\ge 1.
\]
Then \(w\ge0\), \(w(r)\to0\) as \(r\to\infty\), and
\[
w'(r)
=-\e^{-br}\left(\frac b r+\frac1{r^2}\right)\le0.
\]
In particular,
\[
\int_1^\infty|w'(r)|\dd r=w(1).
\]
Integration by parts yields
\[
\int_1^\infty w(r)\sin r\dd r
=w(1)\cos 1+\int_1^\infty w'(r)\cos r\dd r.
\]
Therefore,
\[
\left|\int_1^\infty w(r)\sin r\dd r\right|
\le w(1)+\int_1^\infty|w'(r)|\dd r\leq2w(1)
\lesssim 1.
\]
Combining the two parts, we obtain for any $\epsilon>0\text{ and }a\in\mathbb R$ that
\[
|F_\epsilon(a)|\lesssim 1.
\]
Finally, using \eqref{F}, Plancherel's formula, Lemma \ref{l1} the unitarity of \(\mathcal V\) gives
\[
\begin{aligned}
\|M_\epsilon h\|_{L^2}^2
&=\int_{\mathbb R}|F_\epsilon(a)|^2
\bigl(|\mathcal F_z(J_+ h)(a)|^2+|\mathcal F_z(J_- h)(a)|^2\bigr)\dd a\\
&\lesssim\int_{\mathbb R}
\bigl(|\mathcal F_z(J_+ h)(a)|^2+|\mathcal F_z(J_- h)(a)|^2\bigr)\dd a\\
&\lesssim\|h\|_{L^2}^2.
\end{aligned}
\]
Thus, we obtain \eqref{eq:dilation-M-epsilon-uniform}.

For $a\in\mathbb R$, we define
\[
G(a)=\int_0^\infty \e^{-4s}\frac{\sin(as)}s\dd s.
\]
This integral is absolutely convergent for each fixed $a$, since
\[
\left|\e^{-4s}\frac{\sin(as)}s\right|
\le |a|\e^{-4s}.
\]
Consequently, we get that
\[
F_\epsilon(a)\longrightarrow G(a)
\qquad \mbox{as}\ \ \ \ \ \epsilon\to 0.
\]
Notice that
\[
G'(a)
=\int_0^\infty \e^{-4s}\cos(as)\dd s
=\mathfrak{R}\frac1{4-ia}
=\frac4{16+a^2}.
\]
Recall the definition of $F$ in \eqref{eq:dilation-F-M-m}, we get that
\[
G(a)=\int_0^a\frac4{16+\tau^2}\dd \tau
=\arctan\frac{a}{4}=F(a).
\]
Hence for every $a\in\mathbb R$,
\begin{equation}\label{eq:dilation-F-pointwise}
F_\epsilon(a)\longrightarrow F(a)
\qquad \mbox{as}\ \ \ \ \ \epsilon\to 0.
\end{equation}

For $h\in L^2(\mathbb R_\eta)$,  using \eqref{F}--\eqref{eq:dilation-F-pointwise}, \eqref{eq:dilation-F-M-m}, $|\mathcal F_z(J_+ h)|^2+|\mathcal F_z(J_- h)|^2\in L^1(\mathbb R)$ and the unitarity of $\mathcal V$, we obtain that as $\epsilon\to 0$,
\[
\|(M_\epsilon-M)h\|_{L^2}^2
=
\int_{\mathbb R}
|F_\epsilon(a)-F(a)|^2
\bigl(|\mathcal F_z(J_+ h)(a)|^2+|\mathcal F_z(J_- h)(a)|^2\bigr)\dd a \to 0.
\]
This gives
\eqref{eq:dilation-M-strong}.

Finally,  we turn to derive the bounds for the positive multiplier $m$.
Thanks to \eqref{eq:dilation-F-M-m}, we can rewrite $m(a)$ as
\[
m(a)=\frac12\log\left(1+\frac{16}{a^2+9}\right).
\]
For $x\ge0$, there holds
\[
\frac{x}{1+x}\le\log(1+x)\le x.
\]
Applying this inequalities with $x=16/(a^2+9)$ gives
\[
\frac8{a^2+25}
\le m(a)
\le\frac8{a^2+9}.
\]
Finally, notice that
\[
\frac8{a^2+25}\ge\frac4{a^2+16},
\qquad
\frac8{a^2+9}\le\frac{16}{a^2+16},
\]
we derive \eqref{eq:dilation-m-bounds}. This completes the proof this lemma.
\end{proof}

The following lemma establishes the compatibility of the bounded multiplier \(M\) with the weighted norms and derivative norms used below.

\begin{lemma}\label{l3.3}
Let
$M$ and \(M_\epsilon\) be defined as in \eqref{eq:dilation-F-M-m} and \eqref{eq:dilation-M-epsilon}, respectively,
then it holds
\begin{align}
&\|\eta Mh\|_{L^2}\lesssim\|\eta h\|_{L^2}
\qquad\text{if } h\in D(\eta),\label{2.26d}\\
 & \|\partial_\eta(Mh)\|_{L^2}\lesssim\|\partial_\eta h\|_{L^2}
\qquad \text{if } h\in H^1(\mathbb R_\eta),\label{2.27d}  
\end{align}
and 
\begin{align}
  &M_\epsilon h\longrightarrow Mh
\quad\text{in }H^1(\mathbb R_\eta),\ \ \ \ \text{as } \epsilon \to 0,
\qquad\text{if } h\in H^1(\mathbb R_\eta),\label{2.28d}\\
&M_\epsilon h\longrightarrow Mh
\quad\text{in }D(\eta),\ \ \ \ \ \ \   \text{as } \epsilon \to 0,
\qquad\text{if } h\in D(\eta). \label{2.29d} 
\end{align}

Moreover, for any \(\kappa\geq 1\), it holds that
\begin{align}
&\|(\kappa^2+\eta^2)^{1/2}Mh\|_{L^2}
\lesssim\|(\kappa^2+\eta^2)^{1/2}h\|_{L^2}
\qquad
\text{if }(\kappa^2+\eta^2)^{1/2}h\in L^2,\label{2.30d}\\
&\|(\kappa^2+\eta^2)^{-1/2}Mh\|_{L^2}
\lesssim\|(\kappa^2+\eta^2)^{-1/2}h\|_{L^2}
\qquad
\text{if}\ \ h\in L^2.\label{2.31d}
\end{align}
\end{lemma}
\begin{proof}
We first work with \(M_\epsilon\), for which the dilation integral is truncated away from zero. 
For \(s>0\), thanks to Lemma \ref{l1}, a direct calculation yields 
\[
\|\eta\mathcal U_sh-\eta\mathcal U_{-s}h\|_{L^2}
\lesssim(\e^s+\e^{-s})\|\eta h\|_{L^2}.
\]
Using this observation and \eqref{eq:dilation-M-epsilon}, we obtain
\begin{align}
  \eta M_\epsilon h-M_\epsilon(\eta h)=\int_\epsilon^\infty
\frac{\e^{-4s}}{2is}
\left[
(\e^{-s}-1)\mathcal U_s(\eta h)
-(\e^s-1)\mathcal U_{-s}(\eta h)
\right]\dd s .\notag
\end{align}
Notice that for \(s>0\),
\begin{align}
&\int_0^\infty\frac{\e^{-4s}}{2s}
\left(|\e^{-s}-1|+|\e^s-1|\right)\dd s=\int_0^\infty\frac{\e^{-4s}}{2s}(\e^s-\e^{-s})\dd s\notag\\
&=\int_0^\infty\frac{\e^{-3s}-\e^{-5s}}{2s}\dd s\lesssim 1.\notag
\end{align}
Therefore, using Lemma \ref{l1}, we get
\begin{align}\label{2.26c}
\|\eta M_\epsilon h-M_\epsilon(\eta h)\|_{L^2_\eta}&\le \int_\epsilon^\infty
\frac{\e^{-4s}}{2s}
\left\|
(\e^{-s}-1)\mathcal U_s(\eta h)
-(\e^s-1)\mathcal U_{-s}(\eta h)
\right\|_{L^2_\eta}\dd s\notag\\ 
&\le\int_\epsilon^\infty
\frac{\e^{-4s}}{2s}
\left[
|\e^{-s}-1|\|\mathcal U_s(\eta h)\|_{L^2_\eta}
+|\e^s-1|\|\mathcal U_{-s}(\eta h)\|_{L^2_\eta}
\right]\dd s\notag\\ 
&\le\left(
\int_\epsilon^\infty
\frac{\e^{-4s}}{2s}
\bigl[(1-\e^{-s})+(\e^s-1)\bigr]\dd s
\right)\|\eta h\|_{L^2_\eta}\lesssim \|\eta h\|_{L^2_\eta}.
\end{align}
Similarly, it can be proved that
\begin{equation}\label{2.27c}
  \|-i\partial_\eta(M_\epsilon h)-M_\epsilon(-i\partial_\eta h)\|_{L^2_\eta}\lesssim \|\partial_\eta h\|_{L^2_\eta}.  
\end{equation}
In addition, we define 
\begin{align}
  &I_{01}=\int_0^\infty
\frac{\e^{-4s}}{2is}
\left[
(\e^{-s}-1)\mathcal U_s(\eta h)
-(\e^s-1)\mathcal U_{-s}(\eta h)
\right]\dd s ,\notag\\
&I_{02}
=\int_0^\infty
\frac{\e^{-4s}}{2is}
\left[
(\e^s-1)\mathcal U_s(-i\partial_\eta h)
-(\e^{-s}-1)\mathcal U_{-s}(-i\partial_\eta h)
\right]\dd s.\notag
\end{align}
Similar to the above, by a straightforward calculation, we can prove
\begin{align}\label{2.27aa}
&\lim\limits_{\epsilon \to 0}\|(\eta M_\epsilon h-M_\epsilon(\eta h))-I_{01}\|_{L^2}=0,\notag\\
&\lim\limits_{\epsilon \to 0}\|(-i\partial_\eta(M_\epsilon h)-M_\epsilon(-i\partial_\eta h))-I_{02}\|_{L^2}=0.
\end{align}
By \eqref{eq:dilation-M-strong} in Lemma \ref{l2}, we know that
\begin{equation}\label{2.28aa}
    \lim_{\epsilon\to0}
\|(M_\epsilon-M)h\|_{L^2}=0,\ \ \ \lim_{\epsilon\to0}
\|(M_\epsilon-M)(\eta h)\|_{L^2}=0.
\end{equation}
Thus, we deduce from \eqref{2.27aa} and \eqref{2.28aa} that 
\begin{align}\label{2.29aa}
   &\lim\limits_{\epsilon \to 0}\|\eta M_\epsilon h-M(\eta h)-I_{01}\|_{L^2} \notag\\ 
   \lesssim& \lim\limits_{\epsilon \to 0}\|(M_\epsilon(\eta h)-M(\eta h))+((\eta M_\epsilon h-M_\epsilon(\eta h))-I_{01})\|_{L^2}\notag\\
   \lesssim& \lim\limits_{\epsilon \to 0}\|M_\epsilon(\eta h)-M(\eta h)\|_{L^2}+\lim\limits_{\epsilon \to 0}\|(\eta M_\epsilon h-M_\epsilon(\eta h))-I_{01}\|_{L^2}=0.
\end{align}
As the multiplication operator $\eta$ is closed, this yields
\begin{equation}\notag
    I_{01}=\eta M h-M(\eta h).
\end{equation}
By a similar argument, we can also prove that
\begin{equation}\notag
    I_{02}=-i\partial_\eta(M h)-M(-i\partial_\eta h).
\end{equation}
Therefore, from \eqref{2.28aa} and \eqref{2.29aa}, we derive
\begin{align}\label{2.30aa}
&\lim\limits_{\epsilon \to 0}\|(\eta M_\epsilon h-M_\epsilon(\eta h))-(\eta M h-M(\eta h))\|_{L^2_\eta}=0,\notag\\
&\lim\limits_{\epsilon \to 0}\|(-i\partial_\eta(M_\epsilon h)-M_\epsilon(-i\partial_\eta h))-(-i\partial_\eta(M h)-M(-i\partial_\eta h))\|_{L^2_\eta}=0,
\end{align}
and then
\begin{align}\label{2.31aa}
&\lim\limits_{\epsilon \to 0}\|(\eta M_\epsilon h-\eta M h\|_{L^2_\eta}=0,\notag\\
&\lim\limits_{\epsilon \to 0}\|\partial_\eta(M_\epsilon h)-\partial_\eta(M h)\|_{L^2_\eta}=0,
\end{align}
which together with \eqref{eq:dilation-M-strong} gives \eqref{2.28d}--\eqref{2.29d}.

Similar to \eqref{2.26c}, we can prove
\begin{equation}\notag
    \|\eta M h-M(\eta h)\|_{L^2}\lesssim \|\eta h\|_{L^2},
\end{equation}
then by \eqref{eq:dilation-M-basic}, we get that
\[
\|\eta Mh\|_{L^2}\lesssim \|M(\eta h)\|_{L^2}+\|(\eta M h-M(\eta h)\|_{L^2}\lesssim \|\eta h\|_{L^2}.
\]
This gives \eqref{2.26d}. We can prove \eqref{2.27d} similarly.

Next, for fixed \(\kappa \geq 1\), the condition \((\kappa^2+\eta^2)^{1/2}h\in L^2\) is equivalent to \(h\in D(\eta)\). Using
\[
\|(\kappa^2+\eta^2)^{1/2}h\|_{L^2}^2
=\kappa^2\|h\|_{L^2}^2+\|\eta h\|_{L^2}^2,
\]
we obtain
\[
\begin{aligned}
\|(\kappa^2+\eta^2)^{1/2}Mh\|_{L^2}^2
&=\kappa^2\|Mh\|_{L^2}^2+\|\eta Mh\|_{L^2}^2\\
&\lesssim\bigl(\kappa^2\|h\|_{L^2}^2+\|\eta h\|_{L^2}^2\bigr)\lesssim\|(\kappa^2+\eta^2)^{1/2}h\|_{L^2}^2.
\end{aligned}
\]
This yields \eqref{2.30d}.

Let \(W:=(\kappa^2+\eta^2)^{1/2}\). For fixed \(\kappa\geq 1\), \(W^{-1}\) is a bounded self-adjoint multiplication operator. Thus, for \(h\in L^2\),
\[
\|W^{-1}Mh\|_{L^2}
=
\sup_{\substack{f\in D(W)\\\|Wf\|_{L^2}=1}}
|\langle Mh,f\rangle|.
\]
By self-adjointness of \(M\) and H\"older's inequality,
\[
\begin{aligned}
|\langle Mh,f\rangle|
&=|\langle h,Mf\rangle|=|\langle W^{-1}h,WMf\rangle|\\
&\le\|W^{-1}h\|_{L^2}\,\|WMf\|_{L^2}\lesssim\|W^{-1}h\|_{L^2}\,\|Wf\|_{L^2}.
\end{aligned}
\]
Taking the supremum over \(\|Wf\|_{L^2}=1\) yields
\[
\|(\kappa^2+\eta^2)^{-1/2}Mh\|_{L^2}
\le2\|(\kappa^2+\eta^2)^{-1/2}h\|_{L^2}.
\]
Thus, we get \eqref{2.31d}, and complete the proof of this lemma.
\end{proof}

After introducing the multiplier \(M\), we next investigate the properties of the commutator of \(M\) with the operator \(i \left(-\partial_\eta^2+\frac{2}{|k|^2 + \eta^2}\right)\). 
Thus, for any $f,g\in H^1(\mathbb R)$, we define
\begin{equation}\label{b}
\langle i[H_\kappa,M]f,g\rangle
:=\mathfrak b_\kappa(f,g)=i\bigl(\mathfrak h_\kappa(Mf,g)-\mathfrak h_\kappa(f,Mg)\bigr),
\end{equation}
where $p:=-i\partial_\eta$, $\kappa\geq 1$ and 
\[
H_\kappa=p^2+\frac{2}{\kappa^2+\eta^2},\ \ \ \ \mathfrak h_\kappa(f,g)
=\langle pf,pg\rangle
+\left\langle \frac{2}{\kappa^2+\eta^2}f,g\right\rangle.
\]
This definition agrees with the usual operator commutator, by the self-adjointness of \(M\).
Moreover, the Hermitian symmetry of \(\mathfrak h_\kappa\) implies
\[
\mathfrak b_\kappa(f,g)
=\overline{\mathfrak b_\kappa(g,f)}.
\]

We now begin to analyze the properties of this commutator.

\begin{lemma}\label{l3.4}
  Let $M$ and $m(a)$ be defined by \eqref{eq:dilation-F-M-m} and $p=i\partial_\eta$,
then for all \(f,g\in H^1(\mathbb R_\eta)\),
\begin{equation}\label{2.42c}
i\bigl(\langle pMf,pg\rangle-\langle pf,pMg\rangle\bigr)
=\langle m(A)pf,pg\rangle,
\end{equation}
or equivalently,
\begin{equation}\label{2.43c}
i[p^2,M]=p\,m(A)\,p
\end{equation}
as an identity of continuous operators from \(H^1\) to \(H^{-1}\).
In particular, the kinetic commutator is nonnegative:
\begin{equation}\label{2.44c}
i\bigl(\langle pMh,ph\rangle-\langle ph,pMh\rangle\bigr)
=\|m(A)^{1/2}ph\|_{L^2}^2\ge0,
\qquad h\in H^1(\mathbb R_\eta).
\end{equation}
\end{lemma}
\begin{proof}
    
Using \eqref{3.6bb} in Lemma \ref{l1}, we know that
\[
p\mathcal U_s=\e^s\mathcal U_sp
\]
on \(H^1(\mathbb R_\eta)\). On \(H^2(\mathbb R_\eta)\), we have
\[
p^2\mathcal U_s=\e^{2s}\mathcal U_sp^2,
\qquad
p\mathcal U_sp=\e^s\mathcal U_sp^2.
\]
Consequently,
\[
[p^2,\mathcal U_s]
=(\e^{2s}-1)\mathcal U_sp^2=(\e^s-\e^{-s})p\mathcal U_sp=2\sinh s\,p\mathcal U_sp.
\]
Replacing \(s\) by \(-s\), we likewise obtain
\[
[p^2,\mathcal U_{-s}]
=-2\sinh s\,p\mathcal U_{-s}p.
\]
To use these identities on \(H^1\), we formulate them directly as sesquilinear forms. For \(f,g\in H^1\), set
\[
\mathfrak c_s(f,g)
=\langle p\mathcal U_sf,pg\rangle
-\langle pf,p\mathcal U_{-s}g\rangle.
\]
The appearance of \(\mathcal U_{-s}\) in the second term is due to
\(\mathcal U_s^*=\mathcal U_{-s}\). Using the dilation identities and unitarity,
\[
\begin{aligned}
\langle p\mathcal U_sf,pg\rangle
&=\e^s\langle\mathcal U_spf,pg\rangle,\\
\langle pf,p\mathcal U_{-s}g\rangle
&=\e^{-s}\langle pf,\mathcal U_{-s}pg\rangle=\e^{-s}\langle\mathcal U_spf,pg\rangle.
\end{aligned}
\]
Hence
\begin{equation}\label{c1}
\mathfrak c_s(f,g)
=2\sinh s\,\langle\mathcal U_spf,pg\rangle.
\end{equation}
Similarly,
\begin{equation}\label{c2}
\mathfrak c_{-s}(f,g)
=-2\sinh s\,\langle\mathcal U_{-s}pf,pg\rangle.
\end{equation}

For fixed \(\epsilon>0\), we derive from \eqref{eq:dilation-M-epsilon} and \eqref{c1}--\eqref{c2} that
\begin{align}\label{2.41a}
&i\bigl(\langle pM_\epsilon f,pg\rangle
       -\langle pf,pM_\epsilon g\rangle\bigr)\notag\\
       &=\int_\epsilon^\infty
\frac{\e^{-4s}}{2s}
\left(
\langle p\mathcal U_sf,pg\rangle
-\langle p\mathcal U_{-s}f,pg\rangle
\right)\dd s+\int_\epsilon^\infty
\frac{\e^{-4s}}{2s}
\left(
\langle pf,p\mathcal U_sg\rangle
-\langle pf,p\mathcal U_{-s}g\rangle
\right)\dd s  \notag\\
&=
\int_\epsilon^\infty
\frac{\e^{-4s}}{2s}
\bigl(\mathfrak c_s(f,g)-\mathfrak c_{-s}(f,g)\bigr)\dd s\notag\\
&=
\int_\epsilon^\infty
\e^{-4s}\frac{\sinh s}{s}
\langle(\mathcal U_s+\mathcal U_{-s})pf,pg\rangle\dd s.
\end{align}
Define
\[
N_\epsilon u
=\int_\epsilon^\infty
\e^{-4s}\frac{\sinh s}{s}
(\mathcal U_su+\mathcal U_{-s}u)\dd s,
\qquad u\in L^2.
\]
Then
\begin{equation}\label{N1}
i\bigl(\langle pM_\epsilon f,pg\rangle
       -\langle pf,pM_\epsilon g\rangle\bigr)
=\langle N_\epsilon pf,pg\rangle.
\end{equation}
Since $p$ is a self-adjoint operator, we deduce
\[
i[p^2,M_\epsilon]=pN_\epsilon p.
\]
Define 
\[
Nu=\int_0^\infty
\e^{-4s}\frac{\sinh s}{s}
(\mathcal U_su+\mathcal U_{-s}u)\dd s.
\]
By using Lemma \ref{l1}, we infer
\begin{align}
\|(N-N_\epsilon)u\|_{L^2}
\lesssim&
\left(\int_0^\epsilon \e^{-4s}\frac{\sinh s}{s}\dd s\right)(\|\mathcal U_su\|_{L^2}+\|\mathcal U_{-s}u\|_{L^2})\notag\\ 
\lesssim&
\left(\int_0^\epsilon \e^{-4s}\frac{\sinh s}{s}\dd s\right)\|u\|_{L^2}.
\end{align}
Therefore,
\[
\|N-N_\epsilon\|_{L^2\to L^2}
\le
2\int_0^\epsilon \e^{-4s}\frac{\sinh s}{s}\dd s
\longrightarrow0.
\]
By using \eqref{V}--\eqref{V2}, we have
\[
\begin{aligned}
\mathcal VNh
&=\int_0^\infty
\e^{-4s}\frac{\sinh s}{s}
\,\mathcal V(\mathcal U_s+\mathcal U_{-s})h\dd s\\
&=\int_0^\infty
\e^{-4s}\frac{\sinh s}{s}
\bigl(2\cos(as)F_z(J_+ h)(a),\,2\cos(as)F_z(J_- h)(a)\bigr)\dd s\\
&=\bigl(n(a)F_z(J_+ h)(a),\,n(a)F_z(J_- h)(a)\bigr),
\end{aligned}
\]
where
\[
n(a)
=2\int_0^\infty
\e^{-4s}\frac{\sinh s}{s}\cos(as)\dd s.
\]
Consequently, by \eqref{ah}, we have \(N=n(A)\). In addition, 
by using 
\[
2\e^{-4s}\sinh s=\e^{-3s}-\e^{-5s},\ \ \frac{\e^{-3s}-\e^{-5s}}s
=\int_3^5\e^{-\theta s}\dd \theta,
\]
we infer
\[
\begin{aligned}
n(a)
&=\int_3^5
\left(\int_0^\infty \e^{-\theta s}\cos(as)\dd s\right)d\theta=\int_3^5\frac{\theta}{\theta^2+a^2}\dd \theta\\
&=\frac12\log\frac{a^2+25}{a^2+9}=m(a).
\end{aligned}
\]
We have therefore proved
\begin{equation}\label{N}
N=m(A),
\qquad
N_\epsilon\longrightarrow m(A)
\quad\text{in } L^2.
\end{equation}

Let \(f,g\in H^1(\mathbb R_\eta)\), then it follows from \eqref{2.28d} that
\[
\langle pM_\epsilon f,pg\rangle\to\langle pMf,pg\rangle, \ \ \ \ \ \epsilon \to 0,
\]
and
\[
\langle pf,pM_\epsilon g\rangle\to\langle pf,pMg\rangle  \ \ \ \ \ \epsilon \to 0.
\]
On the other hand, thanks to \eqref{N},
\[
\bigl|\langle(N_\epsilon-m(A))pf,pg\rangle\bigr|
\le
\|N_\epsilon-m(A)\|_{L^2\to L^2}\|pf\|_{L^2}\|pg\|_{L^2}
\to0,\  \ \ \epsilon \to 0,
\]
which together with \eqref{N1} shows
\begin{equation}\label{pM}
i\bigl(\langle pMf,pg\rangle-\langle pf,pMg\rangle\bigr)
=\langle m(A)pf,pg\rangle.
\end{equation}
Thus,we obtain \eqref{2.42c}--\eqref{2.43c}.
Finally,
taking \(f=g=h\) in \eqref{2.42c}, we obtain \eqref{2.44c}.
This proves the lemma.

\end{proof}

Recalling the definitions of $F, M$ and $m$ in \eqref{eq:dilation-F-M-m} and of $J$ in \eqref{3.7b}--\eqref{3.8b}, we derive
\[
J_\pm\left(\frac{2}{1+\eta^2}h\right)(z)
=\e^{z/2}\frac{2}{1+(\pm \e^z)^2}h(\pm \e^z)
=\frac{2}{1+\e^{2z}}(J_\pm h)(z),
\]
and 
\[
\widehat{J_\pm(Mh)}(\xi)
=
\arctan(\xi/4)\,
\widehat{J_\pm h}(\xi),
\]
that is,
\[
J(\frac{2}{1+\eta^2})J^{-1}=W\oplus W,
\qquad
JMJ^{-1}=\widetilde M\oplus\widetilde M,
\]
where
\begin{equation}\label{W}
W(z)=\frac{2}{1+\e^{2z}}=1-\tanh z,
\qquad
\widetilde M=F(-i\partial_z).
\end{equation}
Here \(-i\partial_z\) is self-adjoint on \(H^1(\mathbb R_z)\) and \(F(-i\partial_z)\) is bounded and self-adjoint on all of \(L^2(\mathbb R_z)\).

\begin{lemma}\label{l5}
Define 
\[
K(r)=
\begin{cases}
\displaystyle \frac{\sinh r}{2r}\e^{-4|r|},&r\ne0,\\[2mm]
\displaystyle \frac12,&r=0,
\end{cases}
\]
then \(K\in L^1(\mathbb R)\), and the convolution operator
\[
C_Kg=K*g
\]
satisfies
\begin{equation}\label{2.52c}
i[W,F(-i\partial_z)]=\operatorname{sech}z\  C_K \operatorname{sech}z.
\end{equation}
 Moreover,
\begin{equation}\label{2.53c}
\int_{\mathbb R}\e^{-i\xi r}K(r)\dd r
=
\frac14\log\frac{\xi^2+25}{\xi^2+9}
=
\frac12m(\xi), \ \ C_K=\frac12m(-i\partial_z),
\end{equation}
where $m$ is defines in \eqref{eq:dilation-F-M-m}.   Consequently,
 for every \(g\in L^2(\mathbb R_z)\),
\begin{equation}\label{2.54c}
\left\langle i[W,F(-i\partial_z)]g,g\right\rangle
=
\frac12
\left\|m(-i\partial_z)^{1/2}(\operatorname{sech}z\ g)\right\|_{L^2}^2
\ge0.
\end{equation}

\end{lemma}

\begin{proof}
We first establish the kernel identity for the truncated multiplier. For \(\epsilon>0\) and $g(z)\in L^2$, let
\[
\widetilde M_\epsilon
=
F_\epsilon(-i\partial_z),\qquad
(T_sg)(z)=g(z+s),
\]
where $F_\epsilon$ is defined in \eqref{eq:dilation-F-epsilon}. By using 
\[
\begin{aligned}
\mathcal F_z\left[
\int_\epsilon^\infty
\frac{\e^{-4s}}{2is}(T_sg-T_{-s}g)\dd s
\right](\xi)
&=
\left(
\int_\epsilon^\infty
\e^{-4s}\frac{\sin(s\xi)}s\dd s
\right)\widehat g(\xi)\\
&=F_\epsilon(\xi)\widehat g(\xi),
\end{aligned}
\]
we get that
\[
\widetilde M_\epsilon g
=F_\epsilon(-i\partial_z)g=
\int_\epsilon^\infty
\frac{\e^{-4s}}{2is}
\bigl(T_sg-T_{-s}g\bigr)\dd s.
\]
For each fixed \(\xi\in \mathbb{R}\), note that
\[
\left|\e^{-4s}\frac{\sin(s\xi)}s\right|
\le |\xi|\e^{-4s},
\]
thus, we deduce 
\[
F_\epsilon(\xi)\longrightarrow
\int_0^\infty \e^{-4s}\frac{\sin(s\xi)}s\dd s
=\arctan\frac{\xi}{4}=F(\xi), \ \ \ \epsilon \to 0,
\]
where $F(\xi)$ is defined in \eqref{eq:dilation-F-M-m}. Then by Plancherel's formula and the boundness of $F_\epsilon(\xi)$ and $F(\xi)$, we infer
\[
\begin{aligned}
\|\widetilde M_\epsilon g-F(-i\partial_z)g\|_{L^2}^2
&=
\int_{\mathbb R}
|F_\epsilon(\xi)-F(\xi)|^2
|\widehat g(\xi)|^2\dd \xi \to0,\ \ \ \epsilon \to 0,
\end{aligned}
\]
for every \(g\in L^2\).

In addition, by changing variables separately in the two translation terms, we have
\begin{align}\label{Mg}
(\widetilde M_\epsilon g)(z)
&=
\int_\epsilon^\infty
\frac{\e^{-4s}}{2is}
\bigl(g(z+s)-g(z-s)\bigr)\dd s=
\int_{\mathbb R}
L_\epsilon(z-z')g(z')\dd z',
\end{align}
where
\[
L_\epsilon(r)
=
\frac{i}{2r}\e^{-4|r|}
\mathbf1_{\{|r|\ge\epsilon\}}.
\]
Hence, 
\[
(\widetilde M_\epsilon g)(z)=C_{L_\epsilon}=(L_\epsilon*g)(z).
\]
Using \eqref{Mg}, it holds for \(g\in L^2\) that
\[
\begin{aligned}
\bigl(i[W,\widetilde M_\epsilon]g\bigr)(z)
=
\int_{\mathbb R}
i\bigl(W(z)-W(z')\bigr)
L_\epsilon(z-z')g(z')\dd z'.
\end{aligned}
\]
Since
\[
\begin{aligned}
W(z)-W(z')
&=-\bigl(\tanh z-\tanh z'\bigr)\\
&=-\frac{\sinh z\cosh z'-\cosh z\sinh z'}
        {\cosh z\cosh z'}\\
&=-\frac{\sinh(z-z')}{\cosh z\cosh z'}.
\end{aligned}
\]
By writing \(r=z-z'\) and  
$S(z)=\operatorname{sech}z$, it follows that
\[
\begin{aligned}
i\bigl(W(z)-W(z')\bigr)L_\epsilon(r)
&=
i\left(-\frac{\sinh r}{\cosh z\cosh z'}\right)
\frac{i}{2r}\e^{-4|r|}
\mathbf1_{\{|r|\ge\epsilon\}}\\
&=
S(z)K(r)\mathbf1_{\{|r|\ge\epsilon\}}S(z').
\end{aligned}
\]
Thus, if
\[
K_\epsilon(r)
:=K(r)\mathbf1_{\{|r|\ge\epsilon\}},
\qquad
K(r)=\frac{\sinh r}{2r}\e^{-4|r|},
\]
then it holds
\[
\begin{aligned}
\bigl(i[W,\widetilde M_\epsilon]g\bigr)(z)
&=\int_{\mathbb R}
S(z)K_\epsilon(z-z')S(z')g(z')\dd z'\\
&=S(z)\bigl[K_\epsilon*(Sg)\bigr](z).
\end{aligned}
\]
Moreover, using $
\frac{\sinh r}{r}\to1\quad(r\to0)$ and $|K(r)|
\le
\frac{\e^{-3|r|}}{4|r|}$ for \(|r|\ge1\), we get \(K\in L^1(\mathbb R)\), and then
\[
\begin{aligned}
\|S C_{K_\epsilon}S-S C_KS\|_{L^2\to L^2}
&\le
\|S\|_{L^\infty}^2
\|K_\epsilon-K\|_{L^1}\to0, \ \ \ \ \ \epsilon \to 0.
\end{aligned}
\]
On the other hand, note that \(W\) is bounded  for every \(g\in L^2\),
\[
\begin{aligned}
&\bigl\|
i[W,\widetilde M_\epsilon]g-i[W,F(-i\partial_z)]g
\bigr\|_{L^2}\\
&\qquad\le
\|W\|_{L^\infty}
\|(\widetilde M_\epsilon-F(-i\partial_z))g\|_{L^2}
+
\|(\widetilde M_\epsilon-F(-i\partial_z))(Wg)\|_{L^2}
\to0, \ \ \epsilon \to 0.
\end{aligned}
\]
Uniqueness of the limit yields
\[
i[W,F(-i\partial_z)]=S C_KS=\operatorname{sech}z\  C_K \operatorname{sech}z.
\]
This gives \eqref{2.52c}.

It remains to identify the Fourier multiplier of \(C_K\). 
Using
\[
\frac{\sinh r}{r}
=
\int_0^1\cosh(\theta r)\dd \theta,
\]
 we obtain
\[
K(r)
=
\frac12\int_0^1\e^{-4|r|}\cosh(\theta r)\dd \theta.
\]
 Therefore,
\[
\begin{aligned}
\int_{\mathbb R}\e^{-i\xi r}K(r)\dd r
&=
\frac12\int_0^1
\int_{\mathbb R}
\e^{-i\xi r}\e^{-4|r|}\cosh(\theta r)\dd r\dd \theta.
\end{aligned}
\]
For \(\alpha>0\),
\[
\int_{\mathbb R}\e^{-i\xi r}\e^{-\alpha|r|}\dd r
=
2\int_0^\infty \e^{-\alpha r}\cos(\xi r)\dd r
=
\frac{2\alpha}{\alpha^2+\xi^2}.
\]
Since
\[
\e^{-4|r|}\cosh(\theta r)
=
\frac12
\left(
\e^{-(4-\theta)|r|}
+
\e^{-(4+\theta)|r|}
\right),
\]
we conclude that
\[
\begin{aligned}
\int_{\mathbb R}\e^{-i\xi r}K(r)\dd r
&=
\frac12\int_0^1
\left[
\frac{4-\theta}{(4-\theta)^2+\xi^2}
+
\frac{4+\theta}{(4+\theta)^2+\xi^2}
\right]\dd\theta\\
&=
\frac12\int_3^5\frac{t}{t^2+\xi^2}\dd t\\
&=
\frac14\log\frac{\xi^2+25}{\xi^2+9}=
\frac12m(\xi).
\end{aligned}
\]
 Thus, we get for every \(g\in L^2\) that
\[
\mathcal F_z(C_Kg)(\xi)
=
\left(\int_{\mathbb R}\e^{-i\xi r}K(r)\dd r\right)
(\mathcal F_zg)(\xi)=\frac12m(\xi)(\mathcal F_zg)(\xi).
\]
This means that
\[
C_K=\frac12m(-i\partial_z).
\]
This gives \eqref{2.53c}. In addition, due to \(m\) is bounded, real-valued, and strictly positive on \(\mathbb R\), this is a bounded nonnegative self-adjoint operator. Plancherel’s theorem gives
\[
\langle C_Kg,g\rangle
=
\frac12\int_{\mathbb R}
m(\xi)|(\mathcal F_zg)(\xi)|^2\dd \xi
\ge0,
\]
which together with \eqref{2.52c} gives \eqref{2.54c}.
This completes the proof. 
\end{proof}

\begin{lemma}
   For $A$ as defined in Lemma \ref{l1}, and $M$ and $m$ as defined in \eqref{eq:dilation-F-M-m}--\eqref{eq:dilation-M-basic}, the potential commutator satisfies the bounded-operator identity
\begin{equation}\label{2.56c}
i[\frac{2}{1+\eta^2},M]=2\frac{\eta}{1+\eta^2}\,m(A)\,\frac{\eta}{1+\eta^2}
\end{equation}
on \(L^2(\mathbb R_\eta)\).
Consequently, for every \(f,g\in H^1(\mathbb R_\eta)\),
\begin{equation}\label{2.58c}
\mathfrak b_1(f,g)
=
\langle m(A)\partial_\eta f,\partial_\eta g\rangle
+
2\langle m(A)(\frac{\eta}{1+\eta^2}f),\frac{\eta}{1+\eta^2}g\rangle.
\end{equation}
Here $\mathfrak b_1$ is defined as in \eqref{b}. In particular, for every \(h\in H^1(\mathbb R_\eta)\),
\begin{equation}\label{2.59cc}
\mathfrak b_1(h,h)
=
\|m(A)^{1/2}\partial_\eta h\|_{L^2}^2
+
2\|m(A)^{1/2}(\frac{\eta}{1+\eta^2}h)\|_{L^2}^2
\ge0.
\end{equation}
\end{lemma}
\begin{proof}
 Recalling \eqref{3.7b}--\eqref{3.8b} and \eqref{eq:dilation-F-M-m}--\eqref{eq:dilation-M-basic},
we can first derive
\[
JMJ^{-1}=F(-i\partial_z)\oplus F(-i\partial_z)
\]
and
\[
Jm(A)J^{-1}=m(-i\partial_z)\oplus m(-i\partial_z).
\]
Then by using \eqref{W}--\eqref{2.53c}, we obtain
\[
\begin{aligned}
J\,i[2q^{-1},M]\,J^{-1}
&=
i[W,F(-i\partial_z)]\oplus i[W,F(-i\partial_z)]\\
&=
(\frac12Sm(-i\partial_z) S)\oplus (\frac12Sm(-i\partial_z) S)\\ 
&=\frac12
(S\oplus S)
\bigl(m(D_z)\oplus m(D_z)\bigr)
(S\oplus S).
\end{aligned}
\]
Here we used $S(z)=\operatorname{sech}z$. To identify the multiplication operator corresponding to \(S\), observe that, for \(\eta\ne0\),
\[
S(\log|\eta|)
=
\frac{2}{|\eta|+|\eta|^{-1}}
=
\frac{2|\eta|}{1+\eta^2}.
\]
Hence, defining
\[
\rho(\eta)=\frac{2|\eta|}{1+\eta^2},
\]
we derive
\[
\begin{aligned}
\bigl(J_\pm(\rho h)\bigr)(z)
&=\e^{z/2}\rho(\pm \e^z)h(\pm \e^z)=\rho(\pm \e^z)(J_\pm h)(z)=S(z)(J_\pm h)(z),
\end{aligned}
\]
thus,
\[
J\rho J^{-1}=S\oplus S.
\]
Conjugating the preceding commutator identity by \(J^{-1}\) therefore gives
\begin{align}\label{2.59c}
i[2q^{-1},M]
&=\frac12
\underbrace{J^{-1}(S\oplus S)J}_{\rho}
\underbrace{J^{-1}\bigl(m(D_z)\oplus m(D_z)\bigr)J}_{m(A)}
\underbrace{J^{-1}(S\oplus S)J}_{\rho}\notag\\
&=2(\frac{|\eta|}{1+\eta^2})m(A)(\frac{|\eta|}{1+\eta^2}).
\end{align}

At this stage, the factorization contains \(|a|\). We next show that the absolute values can be removed.

Let \(\Sigma\) denote multiplication by \(\operatorname{sgn}\eta\), that is for $h\in L^2(\eta)$,
\[
(\Sigma h)(\eta)=
\begin{cases}
h(\eta),&\eta>0,\\
-h(\eta),&\eta<0.
\end{cases}
\]
 Its value at \(\eta=0\) is immaterial for its action on \(L^2\). Then
\[
\Sigma^*=\Sigma,
\qquad
\Sigma^2=I,
\]
and
\[
J\Sigma J^{-1}=I\oplus(-I),
\]
where $I$ is the identity operator. Since \(I\oplus(-I)\) commutes with
\(m(-i\partial_z)\oplus m(-i\partial_z)\), it follows that
\[
\Sigma m(A)=m(A)\Sigma.
\]
This commutation reflects the fact that the dilation group preserves the positive and negative half-axes.
Moreover,
\[
\frac{|\eta|}{1+\eta^2}=\frac{\eta}{1+\eta^2}\Sigma=\Sigma \frac{\eta}{1+\eta^2},
\]
Consequently,
\[
\begin{aligned}
\frac{|\eta|}{1+\eta^2}\,m(A)\,\frac{|\eta|}{1+\eta^2}
&=
(\frac{\eta}{1+\eta^2}\Sigma)\,m(A)\,(\Sigma \frac{\eta}{1+\eta^2})\\
&=
\frac{\eta}{1+\eta^2}\,m(A)\,(\Sigma^2\frac{\eta}{1+\eta^2})=
\frac{\eta}{1+\eta^2}\,m(A)\,\frac{\eta}{1+\eta^2}.
\end{aligned}
\]
Combining this with \eqref{2.59c} proves \eqref{2.56c}.

We now derive the full commutator identity. Thanks to \eqref{b}, we have
\[
\begin{aligned}
\mathfrak b_1(f,g)
={}&
i\bigl(
\langle pMf,pg\rangle
-
\langle pf,pMg\rangle
\bigr)\\
&+
i\bigl(
\langle 2q^{-1}Mf,g\rangle
-
\langle 2q^{-1}f,Mg\rangle
\bigr).
\end{aligned}
\]
For the potential contribution, self-adjointness of \(M\) implies
\[
\langle 2q^{-1}f,Mg\rangle
=
\langle M(2q^{-1}f),g\rangle.
\]
Therefore,
\[
\begin{aligned}
&i\bigl(
\langle 2q^{-1}Mf,g\rangle
-
\langle 2q^{-1}f,Mg\rangle
\bigr)=
\langle i[2q^{-1},M]f,g\rangle\\
&\qquad=
2\langle a\,m(A)(af),g\rangle=
2\langle m(A)(af),ag\rangle,
\end{aligned}
\]
whiche together with \eqref{pM} proves \eqref{2.58c}.

\end{proof}

Since
\[
m(\xi)=\frac12\log\frac{\xi^2+25}{\xi^2+9}
\longrightarrow0
\qquad\text{as }|\xi|\longrightarrow\infty,
\]
the nonnegativity of the commutator form does not by itself provide the required control of the weighted norm
\[
\|(\kappa^2+\eta^2)^{-1/2}h\|_{L^2}.
\]
The next step towards a quantitative coercivity estimate is the following identity on the whole line. It controls a weighted \(L^2\)-norm by a derivative term and a multiplication term, with numerical coefficients independent of \(\kappa\).

\begin{lemma}\label{l3.7}
Fix \(\kappa\geq 1\), then for every \(u=u(\eta)\in H^1\),
\begin{equation}\label{2.61c}
\|(\kappa^2+\eta^2)^{-1/2}u\|_{L^2_\eta}^2
\le
\|\partial_\eta u\|_{L^2_\eta}^2+3\|\frac{\eta}{\kappa^2+\eta^2} u\|_{L^2_\eta}^2.
\end{equation}
\end{lemma}

\begin{proof}
For convenience, we denote
\[
q_\kappa(\eta)=\kappa^2+\eta^2,
\qquad
a_\kappa(\eta)=\frac{\eta}{q_\kappa(\eta)}.
\]
 Differentiating \(a_\kappa\) with respect to $\eta$ yields
\[
a_\kappa'(\eta)
=
\frac{(\kappa^2+\eta^2)-2\eta^2}
     {(\kappa^2+\eta^2)^2}
=
\frac{\kappa^2-\eta^2}{q_\kappa(\eta)^2}.
\]
It follows that
\begin{equation}\label{3.63c}
  a_\kappa'(\eta)+2(a_\kappa(\eta))^2=
\frac{\kappa^2-\eta^2+2\eta^2}
     {q_\kappa(\eta)^2}=
\frac{\kappa^2+\eta^2}{q_\kappa(\eta)^2}=
q_\kappa(\eta)^{-1}.  
\end{equation}
 Since \(a_\kappa\) is real-valued,
\[
\|u'+a_\kappa u\|_{L^2}^2
=
\|u'\|_{L^2}^2
+
\|a_\kappa u\|_{L^2}^2
+
2\operatorname{Re}
\int_{\mathbb R}a_\kappa u'\overline u\dd \eta,
\]
and
\[
2\operatorname{Re}
\int_{\mathbb R}a_\kappa u'\overline u\dd \eta
=
\int_{\mathbb R}a_\kappa\,\partial_\eta|u|^2\dd \eta=
-\int_{\mathbb R}a_\kappa'|u|^2\dd \eta.
\]
We can derive
\[
\begin{aligned}
\|u'+a_\kappa u\|_{L^2}^2
&=
\|u'\|_{L^2}^2
+
\|a_\kappa u\|_{L^2}^2
-
\int_{\mathbb R}
\bigl(q_\kappa^{-1}-2a_\kappa^2\bigr)|u|^2\dd \eta\\
&=
\|u'\|_{L^2}^2
+
3\|a_\kappa u\|_{L^2}^2
-
\|q_\kappa^{-1/2}u\|_{L^2}^2.
\end{aligned}
\]
This yields \eqref{2.61c}. 
\end{proof}

\begin{lemma}
For \(c>0\), \(s\in \mathbb{R}\) and \(h\in L^2(\mathbb R)\), define
\begin{equation}\label{3.64c}
  T_ch
=
\int_0^\infty \e^{-cs}\mathcal U_{-s}h\dd s,  
\end{equation}
where \(\mathcal U_{-s}\) is defined in Lemma \ref{l1}. Then \(T_c\) is a bounded linear operator satisfying
\begin{equation}\label{3.65c}
    \|T_c\|_{L^2\to L^2}\le\frac1c,\ \ \ T_c=(c+iA)^{-1}.
\end{equation}
Moreover, for every \(h\in H^1(\mathbb R)\) and $p=-i\partial_\eta$,
\begin{align}\label{3.66c}
pT_ch=T_{c+1}(ph), \ \ \|pT_ch\|_{L^2}\le\frac1{c+1}\|ph\|_{L^2}.
\end{align}
In particular, if \(h\in H^1(\mathbb R)\) and \(u=T_3h\), then
\begin{equation}\label{3.67c}
h=(3+iA)u=\eta u'+\frac72u
,
\end{equation}
and
\begin{equation}\label{3.68c}
T_4(\frac{\eta}{1+\eta^2}h)
=
\frac{\eta}{1+\eta^2}u-T_4\left(\frac{2\eta}{q^2}u\right),
\qquad
\|\frac{\eta}{1+\eta^2}u\|_{L^2}\le2\|T_4(\frac{\eta}{1+\eta^2}h)\|_{L^2}.
\end{equation}
\end{lemma}
\begin{proof}
Firstly, notice that
\[
\int_0^\infty
\|\e^{-cs}\mathcal U_{-s}h\|_{L^2}\dd s
=
\int_0^\infty \e^{-cs}\|h\|_{L^2}\dd s
=
\frac1c\|h\|_{L^2}.
\]
Thus 
\[
\|T_ch\|_{L^2}\le\frac1c\|h\|_{L^2}.
\]
Next, by  \eqref{2.13b}, we infer
\begin{equation}\label{3.69c}
J_\pm T_ch
=
\int_0^\infty \e^{-cs}(J_\pm h)(\,\cdot-s)\dd s=k_c*(J_\pm h),
\end{equation}
where
\[
k_c(s)=\e^{-cs}\mathbf1_{(0,\infty)}(s).
\]
Let \(\mathcal F_z\) denote the Fourier transform, with convention
\[
(\mathcal F_zf)(a)
=
\frac1{\sqrt{2\pi}}
\int_{\mathbb R}\e^{-iza}f(z)\dd z.
\]
 The Fourier multiplier of convolution with \(k_c\) is
\[
\int_{\mathbb R}\e^{-ia s}k_c(s)\dd s
=
\int_0^\infty \e^{-(c+ia)s}\dd s
=
\frac1{c+ia}.
\]
Consequently,
\begin{equation}\label{3.70c}
\mathcal F_z(k_c*(J_\pm h))(a)
=
\frac1{c+ia}\,(\mathcal F_z(J_\pm h))(a).
\end{equation}
Introduce the unitary map
\[
\mathcal V=(\mathcal F_z\oplus\mathcal F_z)J
\]
and write
\[
\mathcal Vh=(\mathcal F_z(J_+h),\mathcal F_z(J_-h)):=(v_+,v_-).
\]
Combining \eqref{ah} and \eqref{3.69c}--\eqref{3.70c},  we arrive at
\[
\begin{aligned}
\mathcal V\bigl((c+iA)T_ch\bigr)
&=
\left(
(c+ia)\frac{v_+}{c+ia},
(c+ia)\frac{v_-}{c+ia}
\right)\\
&=
(v_+,v_-)=
\mathcal Vh.
\end{aligned}
\]
By injectivity of \(\mathcal V\),
\[
(c+iA)T_ch=h.
\]
This proves \eqref{3.65c}.
 
For \(h\in H^1(\mathbb R_\eta)\), by applying \eqref{3.6bb} and \eqref{3.65c}, 
\[
\begin{aligned}
pT_ch
&=
\int_0^\infty \e^{-cs}p\mathcal U_{-s}h\dd s\\
&=
\int_0^\infty \e^{-(c+1)s}\mathcal U_{-s}(ph)\dd s\\
&=
T_{c+1}(ph),
\end{aligned}
\]
and then,
\[
\|pT_ch\|_{L^2}
\le
\frac1{c+1}\|ph\|_{L^2}.
\]
Therefore, we obtain \eqref{3.66c}. 

Now let \(u=T_3h\), it follows from \eqref{3.6b}, \eqref{3.65c}--\eqref{3.66c} that
\[
iAu=\eta u'+\frac12u,
\]
we obtain
\[
h=\eta u'+\frac72u.
\]
Thus \eqref{3.67c} follows.

It remains to prove the multiplication identity. Denote
\[
a(\eta)=\frac{\eta}{1+\eta^2},
\qquad
a'(\eta)=\frac{1-\eta^2}{(1+\eta^2)^2},\ \ \ \ q(\eta)=1+\eta^2.
\]
It follows that
\[
\begin{aligned}
(4+iA)(au)
&=
4au+a\,iAu+\eta a'u\\
&=
a(3+iA)u+(a+\eta a')u\\
&=
ah+\frac{2\eta}{q^2}u.
\end{aligned}
\]
Applying \eqref{3.65c}
to \(au\), we conclude that
\[
T_4(ah)
=
au-T_4\left(\frac{2\eta}{q^2}u\right).
\]
Finally, since \(q\ge1\),
\[
\left|\frac{2\eta}{q^2}\right|
=
\frac2q|a|
\le2|a|.
\]
Together with \(\|T_4\|_{L^2\to L^2}\le\frac{1}{4}\), this gives
\[
\|au-T_4(ah)\|_{L^2}
\le
\frac14\left\|\frac{2\eta}{q^2}u\right\|_{L^2}\le
\frac12\|au\|_{L^2}.
\]
Hence
\[
\|au\|_{L^2}
\le
\|T_4(ah)\|_{L^2}+\frac12\|au\|_{L^2},
\]
and absorption yields
\[
\|au\|_{L^2}\le2\|T_4(ah)\|_{L^2}.
\]
This proves \eqref{3.68c} and completes the proof.
\end{proof}

\begin{lemma}\label{l3.9}
   For every \(h\in H^1(\mathbb R_\eta)\),
\begin{equation}\label{3.72cc}
  \|(1+\eta^2)^{-1/2}h\|_{L^2}^2\lesssim\mathfrak b_1(h,h), 
\end{equation}
here $\mathfrak b_1$ is given in \eqref{2.59cc}.
\end{lemma}

\begin{proof}
It follows from \eqref{eq:dilation-m-bounds} and \eqref{3.65c} that
\[
m(\lambda)\ge\frac4{16+\lambda^2},
\qquad \lambda\in\mathbb R, \ \ \ \ \ T_4=(4+iA)^{-1}.
\]
Then from \eqref{jh}, \eqref{3.69c}--\eqref{3.70c} and Plancherel's formula, we infer
\begin{align}\label{3.72c}
  \|T_4 h\|_{L^2}^2=&\sum_{\sigma=\pm}\|J_\sigma (T_4 h)\|_{L^2}^2=\sum_{\sigma=\pm}\int_\mathbb{R} \left|\frac{\widehat{J_\sigma h}(a)}{4+ia}\right|^2 \dd a=\sum_{\sigma=\pm}\int_\mathbb{R} \frac{|\widehat{J_\sigma h}(a)|^2}{16+a^2} \dd a\notag\\
  \lesssim&\sum_{\sigma=\pm}\int_\mathbb{R} m(a) |\widehat{J_\sigma h}(a)|^2 \dd a\lesssim\langle m(A)h,h\rangle.
\end{align}
By using \eqref{3.72c} and \eqref{2.59cc}, we derive
\begin{equation}
  \|T_4(\partial_\eta h)\|_{L^2}^2+\|T_4(\frac{\eta}{1+\eta^2}h)\|_{L^2}^2\lesssim  \mathfrak b_1(h,h).
\end{equation}
Let $u=T_3h$, then we use \eqref{3.66c}--\eqref{3.68c} and \eqref{2.61c} to get
\begin{align}\notag
  \|(1+\eta^2)^{-1/2}h\|_{L^2}^2\lesssim&\| (1+\eta^2)^{-1/2}\eta \partial_\eta u\|_{L^2}^2+\|(1+\eta^2)^{-1/2}u\|_{L^2}^2\notag\\ 
  \lesssim&\|\partial_\eta u\|_{L^2}^2+\|\frac{\eta}{1+\eta^2} u\|_{L^2_\eta}^2\notag\\
  \lesssim&\|T_4\partial_\eta h\|_{L^2}^2+\|T_4(\frac{\eta}{1+\eta^2}h)\|_{L^2}^2\notag\\ 
  \lesssim&\mathfrak b_1(h,h).
\end{align}
This proves \eqref{3.72cc}.

\end{proof}

We are now able to establish the desired commutator estimate.

\begin{proposition}\label{p3.1}
  For any \(\kappa\geq 1\) and \(f,g\in H^1(\mathbb R_\eta)\), $\mathfrak b_\kappa(f,g)$ and $m$ are defined in \eqref{b} and \eqref{eq:dilation-F-M-m} respectively, we have
\begin{equation}\label{3.76c}
\mathfrak b_\kappa(f,g)
=
\langle m(A)\partial_\eta f,\partial_\eta g\rangle
+
2\langle m(A)(\frac{\eta}{\kappa^2+\eta^2} f),\frac{\eta}{\kappa^2+\eta^2} g\rangle.
\end{equation}
Moreover, 
\begin{equation}\label{3.77c}
\|(\kappa^2+\eta^2)^{-1/2}f\|_{L^2_\eta}^2\lesssim\mathfrak b_\kappa(f,f).
\end{equation} 
\end{proposition}

\begin{proof}
Define
\[
\mathcal D_\kappa:
L^2(\mathbb R_\eta)\longrightarrow L^2(\mathbb R_\xi),
\qquad
(\mathcal D_\kappa h)(\xi)=\kappa^{1/2}h(\kappa\xi).
\]
A change of variables gives
\[
\|\mathcal D_\kappa h\|_{L^2_\xi}^2
=
\int_{\mathbb R}\kappa|h(\kappa\xi)|^2\dd \xi
=
\|h\|_{L^2_\eta}^2.
\]
Its inverse is
\[
(\mathcal D_\kappa^{-1}v)(\eta)
=
\kappa^{-1/2}v(\eta/\kappa).
\]
Thus \(\mathcal D_\kappa\) is unitary. It also maps \(H^1(\mathbb R_\eta)\) bijectively onto \(H^1(\mathbb R_\xi)\).
To distinguish the coordinate representations, write
\[
A_\eta=A=-i\left(\eta\partial_\eta+\frac12\right),
\qquad
A_\xi=-i\left(\xi\partial_\xi+\frac12\right),
\]
where each operator is equipped with its previously established self-adjoint domain in the corresponding \(L^2\)-space. Likewise, set
\[
p_\eta=-i\partial_\eta,
\qquad
p_\xi=-i\partial_\xi,
\qquad
M_\eta=F(A_\eta)=M,
\qquad
M_\xi=F(A_\xi).
\]
For \(h\in H^1(\mathbb R_\eta)\), it holds
\[
p_\xi(\mathcal D_\kappa h)(\xi)
=
-i\kappa^{3/2}h'(\kappa\xi)=
\kappa\,\mathcal D_\kappa(p_\eta h)(\xi).
\]
Thus
\[
p_\xi\mathcal D_\kappa h
=
\kappa\mathcal D_\kappa p_\eta h,
\]
or equivalently,
\[
\mathcal D_\kappa p_\eta\mathcal D_\kappa^{-1}
=
\kappa^{-1}p_\xi.
\]
We next verify how the dilation generator and its multipliers transform. Denote the dilation groups in the two variables by
\[
(\mathcal U_s^{(\eta)}h)(\eta)=\e^{s/2}h(e^s\eta),
\qquad
(\mathcal U_s^{(\xi)}v)(\xi)=\e^{s/2}v(e^s\xi).
\]
For every \(h\in L^2(\mathbb R_\eta)\),
\[
\begin{aligned}
(\mathcal D_\kappa\mathcal U_s^{(\eta)}h)(\xi)
&=
\kappa^{1/2}\e^{s/2}h(\e^s\kappa\xi)\\
&=
(\mathcal U_s^{(\xi)}\mathcal D_\kappa h)(\xi).
\end{aligned}
\]
Consequently,
\[
\mathcal D_\kappa\mathcal U_s^{(\eta)}\mathcal D_\kappa^{-1}
=
\mathcal U_s^{(\xi)}, \ \ \ s\in\mathbb R.
\]
By uniqueness of the self-adjoint generator,
\[
\mathcal D_\kappa D(A_\eta)=D(A_\xi),
\qquad
\mathcal D_\kappa A_\eta\mathcal D_\kappa^{-1}=A_\xi.
\]
Unitary invariance of the functional calculus now gives
\[
\mathcal D_\kappa M_\eta\mathcal D_\kappa^{-1}=M_\xi
\]
and
\[
\mathcal D_\kappa m(A_\eta)\mathcal D_\kappa^{-1}=m(A_\xi).
\]
Equivalently,
\[
M_\xi\mathcal D_\kappa=\mathcal D_\kappa M_\eta,
\qquad
m(A_\xi)\mathcal D_\kappa
=
\mathcal D_\kappa m(A_\eta).
\]
Let
\[
v=\mathcal D_\kappa f,
\quad
w=\mathcal D_\kappa g,\ \ \  q_\kappa(\xi)=\kappa^2+\xi^2,\ \ \ a_\kappa(\xi)=\frac{\xi}{\kappa^2+\xi^2}.
\]
The derivative and multiplication identities, together with unitarity, imply
\[
\begin{aligned}
\mathfrak h_\kappa(f,g)
&=
\kappa^{-2}\langle p_\xi v,p_\xi w\rangle_\xi
+
2\kappa^{-2}\langle q_1^{-1}v,w\rangle_\xi\\
&=
\kappa^{-2}\mathfrak h_1(v,w),
\end{aligned}
\]
where \(\mathfrak h_1\) on the right-hand side is understood in the \(\xi\)-variable.

Using \eqref{b} and
\[
\mathcal D_\kappa M_\eta f=M_\xi v,
\qquad
\mathcal D_\kappa M_\eta g=M_\xi w,
\]
we obtain
\[
\begin{aligned}
\mathfrak b_\kappa(f,g)
&=
i\bigl(
\mathfrak h_\kappa(M_\eta f,g)
-
\mathfrak h_\kappa(f,M_\eta g)
\bigr)\\
&=
i\kappa^{-2}
\bigl(
\mathfrak h_1(M_\xi v,w)
-
\mathfrak h_1(v,M_\xi w)
\bigr)\\
&=
\kappa^{-2}\mathfrak b_1(v,w).
\end{aligned}
\]
Therefore, we arrive at
\[
\mathfrak b_\kappa(f,g)
=
\kappa^{-2}\bigl(
\langle m(A_\xi)p_\xi v,p_\xi w\rangle_\xi+
2\langle m(A_\xi)(a_1v),a_1w\rangle_\xi
\bigr).
\]
We now transfer each term back to the original variable.
For the derivative term, we use
\[
p_\xi v=\kappa\mathcal D_\kappa(p_\eta f),
\qquad
p_\xi w=\kappa\mathcal D_\kappa(p_\eta g)
\]
to obtain
\[
\begin{aligned}
&\kappa^{-2}
\langle m(A_\xi)p_\xi v,p_\xi w\rangle_\xi=
\left\langle
m(A_\xi)\mathcal D_\kappa(p_\eta f),
\mathcal D_\kappa(p_\eta g)
\right\rangle_\xi\\
&\quad=
\left\langle
\mathcal D_\kappa\bigl(m(A_\eta)p_\eta f\bigr),
\mathcal D_\kappa(p_\eta g)
\right\rangle_\xi\\
&\quad=
\langle m(A_\eta)p_\eta f,p_\eta g\rangle_\eta.
\end{aligned}
\]
For the multiplication term, note that
\[
a_1v=\kappa\mathcal D_\kappa(a_\kappa f),
\qquad
a_1w=\kappa\mathcal D_\kappa(a_\kappa g),
\]
hence
\[
\begin{aligned}
&2\kappa^{-2}
\langle m(A_\xi)(a_1v),a_1w\rangle_\xi=
2\left\langle
m(A_\xi)\mathcal D_\kappa(a_\kappa f),
\mathcal D_\kappa(a_\kappa g)
\right\rangle_\xi\\
&\quad=
2\left\langle
\mathcal D_\kappa\bigl(m(A_\eta)(a_\kappa f)\bigr),
\mathcal D_\kappa(a_\kappa g)
\right\rangle_\xi\\
&\quad=
2\langle m(A_\eta)(a_\kappa f),a_\kappa g\rangle_\eta.
\end{aligned}
\]
Combining these computations yields
\[
\mathfrak b_\kappa(f,g)
=
\langle m(A_\eta)p_\eta f,p_\eta g\rangle_\eta
+
2\langle m(A_\eta)(a_\kappa f),a_\kappa g\rangle_\eta.
\]
 Moreover, since \(p_\eta=-i\partial_\eta\),
\[
\begin{aligned}
\langle m(A_\eta)p_\eta f,p_\eta g\rangle_\eta
&=
\langle -i\,m(A_\eta)\partial_\eta f,
        -i\,\partial_\eta g\rangle_\eta\\
&=
(-i)\overline{(-i)}
\langle m(A_\eta)\partial_\eta f,\partial_\eta g\rangle_\eta\\
&=
\langle m(A_\eta)\partial_\eta f,\partial_\eta g\rangle_\eta.
\end{aligned}
\]
Thus, we get \eqref{3.76c}.

Finally, note that
\[
\begin{aligned}
\bigl(\mathcal D_\kappa(q_\kappa^{-1/2}f)\bigr)(\xi)
&=
\kappa^{1/2}
\bigl(\kappa^2+\kappa^2\xi^2\bigr)^{-1/2}
f(\kappa\xi)\\
&=
\kappa^{-1}(1+\xi^2)^{-1/2}v(\xi)\\
&=
\kappa^{-1}q_1^{-1/2}v(\xi).
\end{aligned}
\]
Therefore, by unitarity,
\[
\|q_\kappa^{-1/2}f\|_{L^2_\eta}^2
=
\kappa^{-2}\|q_1^{-1/2}v\|_{L^2_\xi}^2.
\]
Using the unit-frequency bound and the form scaling identity, we conclude that
\[
\|q_\kappa^{-1/2}f\|_{L^2_\eta}^2=
\kappa^{-2}\|q_1^{-1/2}v\|_{L^2_\xi}^2\lesssim\kappa^{-2}\mathfrak b_1(v,v)\lesssim\mathfrak b_\kappa(f,f).
\]
Thus we obtain \eqref{3.77c}. This completes the proof. 
\end{proof}

Following the above preparations, we can establish the following proposition.
\begin{proposition}\label{p3.2}
Let $f_k$ solves \eqref{2.5}. There exists a suitably small positive number $c_1$, such that for $t\geq 0$ and $k\neq 0$,
\begin{align}\label{3.73}
  &\|\e^{c_1\nu^{\frac{1}{2}}t}f_k\|_{L^\infty L^2} +\nu^{\frac{1}{2}}\|\e^{c_1\nu^{\frac{1}{2}}t}\nabla_kf_k\|_{L^2 L^2}+(\nu|k|)^{\frac{1}{4}}\|\e^{c_1\nu^{\frac{1}{2}}t}f_k\|_{L^2 L^2}\notag\\ 
  &+|k|^{\frac{1}{2}}\|\e^{c_1\nu^{\frac{1}{2}}t}(\Delta_k)^{-\frac{1}{2}}f_k\|_{L^2 L^2} \notag\\ 
  \lesssim&\|f_{{\rm in},k}\|_{L^2}+\nu^{-\frac{1}{2}}\|\e^{c_1\nu^{\frac{1}{2}}t} F\|_{L^2L^2}
+\nu^{\frac{1}{2}}
\|\e^{c_1\nu^{\frac{1}{2}}t}G_1\|_{L^2H_{k,y}^{-1}}+|k|^{\frac{1}{2}}\|\e^{c_1\nu^{\frac{1}{2}}t}\nabla_k G_2\|_{L^2L^2}\notag\\ 
&+\nu^{-\frac{1}{2}}\|\e^{c_1\nu^{\frac{1}{2}}t}G_3\|_{L^2H_{k,y}^{-1}},
\end{align}
  where $\nabla_k=(k,\partial_y)$. 
\end{proposition}
\begin{proof}
Fix a nonzero horizontal frequency \(k\), we denote $\sigma=\operatorname{sgn}k$, and choose a fixed constant
\begin{equation}\notag
0<\delta_0\le\frac18,
\end{equation}
independently of \(\nu\), \(k\) and $t$. Define
\begin{equation}\label{3.74}
\mathcal E_\sigma(\hat{f}_k)=\langle  \hat{f}_k-\delta_0\sigma M\hat{f}_k,\hat{f}_k\rangle,
\end{equation}
where multiplier $M$ is defined in \eqref{eq:dilation-F-M-m}. 
Since \(M\) is bounded and self-adjoint, so 
\begin{equation}\notag
\bigl|\langle M\hat{f}_k,\hat{f}_k\rangle\bigr|
\le
\frac{\pi}{2}\|\hat{f}_k\|_{L^2}^2,
\end{equation}
and therefore
\begin{equation}\label{3.75}
\frac{1}{2}\|\hat{f}_k\|_{L^2}^2
\le
\mathcal E_\sigma(\hat{f}_k)
\le
\frac{3}{2}\|\hat{f}_k\|_{L^2}^2.
\end{equation}
Thus the modified energy functional is uniformly equivalent to the ordinary \(L^2\)-energy.
Through direct computation, we obtain
\begin{align}\label{3.76}
\frac{\dd}{\dd t}\mathcal E_\sigma(\hat{f}_k)
=&
2\operatorname{Re} \left\langle-ik\left[-(\partial_\eta)^2+\frac{2}{|k|^2+\eta^2}\right] \hat{f}_k,\hat{f}_k-\delta_0\sigma M\hat{f}_k\right\rangle\notag\\
&-
2\nu\operatorname{Re}\langle (|k|^2+\eta^2) \hat{f}_k,\hat{f}_k-\delta_0\sigma M\hat{f}_k\rangle\notag\\
&+
2\operatorname{Re}
\langle ik\,\hat F_1
+
i\eta\,\hat F_2+\nu \hat G_1-ik\hat G_2+\hat G_3,\hat{f}_k-\delta_0\sigma M\hat{f}_k\rangle.
\end{align}
Recalling \eqref{b}, we obtain
\begin{equation}\notag
\mathfrak h_\kappa(\hat{f}_k,M\hat{f}_k)
=
\langle \partial_\eta \hat{f}_k,\partial_\eta (M\hat{f}_k)\rangle
+
\left\langle\frac{2}{|k|^2+\eta^2}\hat{f}_k,M\hat{f}_k\right\rangle,
\end{equation}
and then 
\begin{equation}\notag
\mathfrak b_\kappa(\hat{f}_k,\hat{f}_k)
=
i\bigl(
\mathfrak h_\kappa(M\hat{f}_k,\hat{f}_k)
-
\mathfrak h_\kappa(\hat{f}_k,M\hat{f}_k)
\bigr)=2\operatorname{Im}\mathfrak h_\kappa(\hat{f}_k,M\hat{f}_k).
\end{equation}
Using \(k\sigma=|k|\), we obtain
\begin{align}\label{3.77}
&2\operatorname{Re} \left\langle-ik\left[-(\partial_\eta)^2+\frac{2}{|k|^2+\eta^2}\right] \hat{f}_k,\hat{f}_k-\delta_0\sigma M\hat{f}_k\right\rangle\notag\\
=&2\operatorname{Re} \left\langle-ik\partial_\eta \hat{f}_k,\partial_\eta(\hat{f}_k-\delta_0\sigma M\hat{f}_k)\right\rangle+2\operatorname{Re} \left\langle-ik\frac{2}{|k|^2+\eta^2} \hat{f}_k,\hat{f}_k-\delta_0\sigma M\hat{f}_k\right\rangle\notag\\
=&2\delta_0|k|\operatorname{Re} \left\langle i\partial_\eta \hat{f}_k,\partial_\eta( M\hat{f}_k)\right\rangle+2\delta_0|k|\operatorname{Re} \left\langle i\frac{2}{|k|^2+\eta^2} \hat{f}_k, M\hat{f}_k\right\rangle\notag\\
=&-2\delta_0|k|\operatorname{Im}\mathfrak h_\kappa(\hat{f}_k,M\hat{f}_k)=-\delta_0|k|\mathfrak b_\kappa(\hat{f}_k,\hat{f}_k).
\end{align}
Thus the choice of sign in \(P_\sigma\) produces the same nonnegative commutator contribution for positive and negative horizontal frequencies.

For the viscous term,
\begin{align}
&-2\nu\operatorname{Re}\langle (|k|^2+\eta^2) \hat{f}_k,\hat{f}_k-\delta_0\sigma M\hat{f}_k\rangle\notag\\
=&
-2\nu\|(|k|^2+\eta^2)^{\frac{1}{2}}h\|_{L^2}^2
+
2\delta_0\sigma\nu
\operatorname{Re}\langle (|k|^2+\eta^2) h,Mh\rangle\notag\\
\leq&-2\nu\|(|k|^2+\eta^2)^{\frac{1}{2}}\hat{f}_k\|_{L^2}^2+2\delta_0\nu\|(|k|^2+\eta^2)^{\frac{1}{2}}\hat{f}_k\|_{L^2}\|(|k|^2+\eta^2)^{\frac{1}{2}}M\hat{f}_k\|_{L^2}\notag\\
\leq&-2\nu\|(|k|^2+\eta^2)^{\frac{1}{2}}\hat{f}_k\|_{L^2}^2+4\delta_0\nu\|(|k|^2+\eta^2)^{\frac{1}{2}}\hat{f}_k\|_{L^2}^2\notag\\
\leq&-\frac{3}{2}\nu\|(|k|^2+\eta^2)^{\frac{1}{2}}\hat{f}_k\|_{L^2}^2.
\end{align}
By \eqref{3.77c}, there exists a constant $0<c_0<1$  such that
\begin{equation}\label{3.79}
2c_0\|(|k|^2+\eta^2)^{-\frac{1}{2}}\hat{f}_k\|_{L^2}^2\leq \mathfrak b_\kappa(\hat{f}_k,\hat{f}_k).
\end{equation}
Notice that
\begin{align}\label{3.80}
   \nu^{\frac{1}{2}}|k|^{\frac{1}{2}} \|\hat{f}_k\|_{L^2}^2
&=\int_{\mathbb R}
\bigl(\nu^{\frac{1}{2}}(|k|^2+\eta^2)^{\frac{1}{2}}|\hat{f}_k|\bigr)
\bigl(|k|^{\frac{1}{2}}(|k|^2+\eta^2)^{-\frac{1}{2}}|\hat{f}_k|\bigr)\dd \eta \notag\\
&\le (\nu^{\frac{1}{2}}\|(|k|^2+\eta^2)^{\frac{1}{2}}\hat{f}_k\|_{L^2})(|k|^{\frac{1}{2}}\|(|k|^2+\eta^2)^{-\frac{1}{2}}\hat{f}_k\|_{L^2})\notag\\ 
&\le \nu\|(|k|^2+\eta^2)^{\frac{1}{2}}\hat{f}_k\|_{L^2}^2+|k|\|(|k|^2+\eta^2)^{-\frac{1}{2}}\hat{f}_k\|_{L^2}^2.
\end{align}
Then it follows that
\begin{align}\label{3.81}
&\frac{\dd}{\dd t}\mathcal E_\sigma(\hat{f}_k)
+\nu\|(|k|^2+\eta^2)^{\frac{1}{2}}\hat{f}_k\|_{L^2}^2+\delta_0 c_0|k|\|(|k|^2+\eta^2)^{-\frac{1}{2}}\hat{f}_k\|_{L^2}^2+\delta_0c_0 \nu^{\frac{1}{2}}|k|^{\frac{1}{2}} \|\hat{f}_k\|_{L^2}^2\notag\\ 
\le&
2\operatorname{Re}
\langle ik\,\hat F_1
+
i\eta\,\hat F_2+\nu \hat G_1-ik\hat G_2+\hat G_3,\hat{f}_k-\delta_0\sigma M\hat{f}_k\rangle.
\end{align}

We next estimate the forcing terms. By using \eqref{2.30d}--\eqref{2.31d}, we infer
\begin{equation}\notag
\|(|k|^2+\eta^2)^{\frac{1}{2}}(\hat{f}_k-\delta_0\sigma M\hat{f}_k)\|_{L^2}
\lesssim\|(|k|^2+\eta^2)^{\frac{1}{2}}\hat{f}_k\|_{L^2}
\end{equation}
and
\begin{equation}\notag
\|(|k|^2+\eta^2)^{-\frac{1}{2}}(\hat{f}_k-\delta_0\sigma M\hat{f}_k)\|_{L^2}
\lesssim\|(|k|^2+\eta^2)^{-\frac{1}{2}}\hat{f}_k\|_{L^2}.
\end{equation}
For the divergence-form forcing, the pointwise Cauchy–Schwarz inequality gives
\begin{equation}\notag
|k\hat F_1(\eta)+\eta\hat F_2(\eta)|^2
\le
(|k|^2+\eta^2)
\bigl(
|\hat F_1(\eta)|^2+|\hat F_2(\eta)|^2
\bigr).
\end{equation}
Thus, by Plancherel’s theorem,
\begin{equation}\notag
\|(|k|^2+\eta^2)^{-\frac{1}{2}}|k\,\hat F_1
+
\eta\,\hat F_2|\|_{L^2}\le\|\hat F\|_{L^2}.
\end{equation}
Consequently,
\begin{align}\label{3.83}
|\langle ik\,\hat F_1
+
i\eta\,\hat F_2,\hat{f}_k-\delta_0\sigma M\hat{f}_k\rangle|
=&
\bigl|
\langle (|k|^2+\eta^2)^{-\frac{1}{2}}(ik\,\hat F_1
+
i\eta\,\hat F_2),
(|k|^2+\eta^2)^{\frac{1}{2}}(\hat{f}_k-\delta_0\sigma M\hat{f}_k)\rangle
\bigr|\notag\\ 
\lesssim&\|\hat F\|_{L^2}\|(|k|^2+\eta^2)^{\frac{1}{2}}\hat{f}_k\|_{L^2}\notag\\
\leq&\varepsilon_0\nu\|(|k|^2+\eta^2)^{\frac{1}{2}}\hat{f}_k\|_{L^2}^2+C_{\varepsilon_0}\nu^{-1}\|\hat F\|_{L^2}^2.
\end{align}
For the forcing terms, a direct calculation yields
\begin{align}\label{3.84}
|\langle\nu \hat G_1,\hat{f}_k-\delta_0\sigma M\hat{f}_k\rangle|
&=
\nu\bigl|
\langle (|k|^2+\eta^2)^{-\frac{1}{2}}\hat g,
(|k|^2+\eta^2)^{\frac{1}{2}}(\hat{f}_k-\delta_0\sigma M\hat{f}_k)\rangle
\bigr|\notag\\
&\lesssim
\nu\|(|k|^2+\eta^2)^{-\frac{1}{2}}\hat G_1\|_{L^2} \|(|k|^2+\eta^2)^{\frac{1}{2}}\hat{f}_k\|_{L^2}\notag\\
&\leq \varepsilon_0\nu\|\nabla_k\hat{f}_k\|_{L^2}^2+C_{\varepsilon_0}\nu
\|G_1\|_{H_{k,y}^{-1}}^2,
\end{align}
\begin{align}\label{3.84a}
|\langle ik\hat G_2,\hat{f}_k-\delta_0\sigma M\hat{f}_k\rangle|
\lesssim& 
|k|\|(|k|^2+\eta^2)^{\frac{1}{2}}\hat{G}_2\|_{L^2}\|(|k|^2+\eta^2)^{-\frac{1}{2}}\hat{f}_k\|_{L^2}\notag\\
\leq& \varepsilon_0 |k|\|(|k|^2+\eta^2)^{-\frac{1}{2}}\hat{f}_k\|_{L^2}^2+ C_{\varepsilon_0}|k|\|\nabla_k G_2\|_{L^2}^2,
\end{align}
and
\begin{align}\label{3.86j}
|\langle\hat G_3,\hat{f}_k-\delta_0\sigma M\hat{f}_k\rangle|
\lesssim&\|(|k|^2+\eta^2)^{-\frac{1}{2}}\hat{G}_3\|_{L^2}\|(|k|^2+\eta^2)^{\frac{1}{2}}\hat{f}_k\|_{L^2}\notag\\ 
\leq& \varepsilon_0 \nu\|(|k|^2+\eta^2)^{\frac{1}{2}}\hat{f}_k\|_{L^2}^2+ C_{\varepsilon_0}\nu^{-1}\|G_3\|_{H_{k,y}^{-1}}^2.
\end{align}
Substituting \eqref{3.83}--\eqref{3.86j} into \eqref{3.81}, choosing $\varepsilon_0$ sufficiently small and using Plancherel's formula, we obtain 
\begin{align}\label{3.86}
&\frac{\dd}{\dd t}\mathcal E_\sigma(f_k)
+\frac{1}{2}\left[\nu\|\nabla_kf_k\|_{L^2}^2+\delta_0 c_0|k|\|(\Delta_k)^{-\frac{1}{2}}f_k\|_{L^2}^2+\delta_0c_0 \nu^{\frac{1}{2}}|k|^{\frac{1}{2}} \|f_k\|_{L^2}^2\right]\notag\\ 
\lesssim&\nu^{-1}\|F\|_{L^2}^2+\nu
\|G_1\|_{H_{k,y}^{-1}}^2+|k|\|\nabla_k G_2\|_{L^2}^2+\nu^{-1}\|G_3\|_{H_{k,y}^{-1}}^2.
\end{align}
Taking $c_1:=\min\{\frac{1}{2},\frac{\delta_0 c_0}{2}\}$, multiplying \eqref{3.86} by $\e^{c_1\nu^{\frac{1}{2}}t}$ and then integrating with respect to $t$, we get \eqref{3.73}, and complete the proof of this proposition.
\end{proof}

\section{Linear estimates for the time-dependent Shear flow \(y^2+V^L(t,y)\)}\label{s4}

In this section, we will consider the space-time estimates for the following linear system:
\begin{equation}\label{4.1}
\begin{cases}
\partial_tf_{k} + \nu \Delta_k f_{k} + i k(V + \partial_y^2V\Delta_k^{-1}) f_{k} = \mathrm{div}_k F, \quad\\[2mm] 
f_{k}(0,y) = f_{k, \mathrm{in}},\ \ \ k\in \mathbb Z~ \mbox{and}~k\neq 0,~y\in \mathbb{R},
\end{cases}
\end{equation}
where \(\Delta_k = -\partial_y^2 + k^2\), \(\mathrm{div}_k F = i k F_1 + \partial_y F_2\).

We can establish the following proposition.
 \begin{proposition}\label{p4.1}
 Let $f_{k}$ solves \eqref{4.1}.   Under the assumptions of Theorem \ref{thm1}, for all $t\geq 0$, we obtain
 \begin{align}\label{4.2}
  &\|\e^{\delta\nu^{\frac{1}{2}}t}f_k\|_{L^\infty L^2} +\nu^{\frac{1}{2}}\|\e^{\delta\nu^{\frac{1}{2}}t}\nabla_kf_k\|_{L^2 L^2}+(\nu|k|)^{\frac{1}{4}}\|\e^{\delta\nu^{\frac{1}{2}}t}f_k\|_{L^2 L^2}\notag\\ 
  &+|k|^{\frac{1}{2}}\|\e^{\delta\nu^{\frac{1}{2}}t}(\Delta_k)^{-\frac{1}{2}}f_k\|_{L^2 L^2} \notag\\ 
  \lesssim&\|f_{k, \mathrm{in}}\|_{L^2}+\nu^{-\frac{1}{2}}\|\e^{\delta\nu^{\frac{1}{2}}t} F\|_{L^2L^2}.
\end{align}
 \end{proposition}

Before proving this proposition, we first present some preliminary results.
For the sake of simplicity, Let's set
\begin{equation}\label{4.1c}
c(t,y)=V_{yy}(t,y)=2+V^L_{yy}(t,y).
\end{equation}
Heat-semigroup contraction and one-dimensional Sobolev embedding give
\begin{equation}\label{4.1cc}
\sup_{t\ge0}\|V^L(t)\|_{W^{5,\infty}}
\le
C\|V_{\mathrm{in}}\|_{H^6}\leq Cc_*\ll 1.
\end{equation}

Since
\[
V_y(t,y)=2y+V^L_y(t,y),
\qquad
V_{yy}(t,y)\ge1,
\]
then there is exactly one critical point \(y_c(t)\), determined by
\[
V_y(t,y_c(t))=0.
\]
It satisfies
\[
|y_c(t)|\le\frac12\|V^L_y(t)\|_{L^\infty}.
\]
Define
\[
\beta(t)=V(t,y_c(t))
\]
and
\begin{equation}\label{4.2c}
Y=\Phi(t,y)
=
\operatorname{sgn}(y-y_c(t))
\sqrt{V(t,y)-\beta(t)}.
\end{equation}
\begin{lemma}\label{l4.1}
For sufficiently small \(c_*\), the map
\(\Phi(t,\cdot):\mathbb R\to\mathbb R\) is an increasing global diffeomorphism. Let
\[
y=\Psi(t,Y),
\qquad
\rho(t,Y)=\Phi_y(t,\Psi(t,Y)),
\qquad
s(t,Y)=\sqrt{c(t,\Psi(t,Y))}.
\]
Then it holds that
\[
V(t,\Psi(t,Y))=Y^2+\beta(t)
\]
and
\begin{equation}\label{4.3cc}
\|\rho-1\|_{W_Y^{3,\infty}}
+
\|s-\sqrt2\|_{W_Y^{3,\infty}}
+
\|s^{-1}-2^{-1/2}\|_{W_Y^{3,\infty}}
\le Cc_*.
\end{equation}
\end{lemma}
\begin{proof}
We denote \(r=y-y_c(t)\), and use Taylor’s formula to get
\[
V(t,y)-V(t,y_c)
=
r^2a_0(t,y),
\]
where
\[
a_0(t,y)
=
\int_0^1(1-\theta)c(t,y_c+\theta r)\dd \theta.
\]
Then by \eqref{4.1c} and \eqref{4.2c}, there holds
\[
a_0=2\int_0^1(1-\theta)\dd \theta+\int_0^1(1-\theta)V^L_{yy}\dd \theta=1+\int_0^1(1-\theta)V^L_{yy}(t,y_c+\theta r)\dd \theta,
\]
\[
\Phi(t,y)=\operatorname{sgn}(r)\sqrt{r^2a_0(t,y)}=r\sqrt{a_0(t,y)}.
\]
For \(0\le j\le3\),
\[
\partial_y^j(a_0-1)
=
\int_0^1
(1-\theta)\theta^j
(V^L)^{(j+2)}(t,y_c+\theta r)\dd \theta,
\]
and then by \eqref{4.1cc},
\[
\|a_0-1\|_{W_y^{3,\infty}}\leq C\|V^L_{yy}\|_{W_y^{3,\infty}}\leq Cc_*.
\]
To estimate \(\Phi_y\) without introducing a factor growing with \(r\), use
\[
\begin{aligned}
V_y(t,y)
&=V_y(t,y)-V_y(t,y_c)=r\int_0^1 V_{yy}(t,y_c+\theta r)\dd \theta\\
&=r\left(2+\int_0^1 V^L_{yy}(t,y_c+\theta r)\dd \theta\right),
\end{aligned}
\]
Since \(\Phi^2=V(t,y)-V(t,y_c)\) and \(2\Phi\Phi_y=V_y\), we obtain
\[
\Phi_y
=
\frac{
1+\frac12\int_0^1V^L_{yy}(t,y_c+\theta r)\dd \theta
}{
\sqrt{a_0(t,y)}
},
\]
and then,
\[
\|\Phi_y-1\|_{W_y^{3,\infty}}=\left\|\frac{
1-\sqrt{a_0(t,y)}+\frac12\int_0^1V^L_{yy}(t,y_c+\theta r)\dd \theta
}{
\sqrt{a_0(t,y)}}\right\|_{W_y^{3,\infty}}\le Cc_*.
\]
Since \(0<c_* \ll 1\), there holds
\[
\frac12\le\Phi_y\le2.
\]
Also, \(\Phi=r\sqrt{a_0}\) tends to \(\pm\infty\) as \(y\to\pm\infty\). Hence \(\Phi\) is a global increasing diffeomorphism. In addition,
since \(\Phi(t,\Psi(t,Y))=Y\), we have for any \(f\in W_y^{3,\infty}\),
\[
\|f\circ\Psi\|_{W_Y^{3,\infty}}
\lesssim \|f\|_{W_y^{3,\infty}}.
\]
Then using \eqref{4.1c}--\eqref{4.1cc} and
\[
\rho-1=(\Phi_y-1)\circ\Psi,\quad
s-\sqrt2=(\sqrt c-\sqrt2)\circ\Psi,\quad
s^{-1}-2^{-1/2}=(c^{-1/2}-2^{-1/2})\circ\Psi,
\]
we get \eqref{4.3cc}, and finish the proof of this lemma.
\end{proof}

Next, we proceed to analyze the perturbed equation.
Denote
\begin{equation}\label{4.3c}
J(t,Y)=\partial_Y\Psi(t,Y)=\rho(t,Y)^{-1},
 \quad \vartheta(t)=\int_0^t\beta(\tau)\dd \tau,
\end{equation}
and introduce the unitary mapping \(
\mathscr U_t:L_y^2\to L_Y^2
\) as
\[
(\mathscr U_t h)(Y)=J(t,Y)^{1/2}h(\Psi(t,Y)).
\]
Then we define
\begin{equation}\label{4.4c}
g_k=\mathcal T_{k,t}f_k,
\qquad
\mathcal T_{k,t}
=\e^{ik\vartheta(t)}\mathscr U_t c^{-1/2}
=\e^{ik\vartheta(t)}s^{-1}\mathscr U_t.
\end{equation}
The last equality follows from \(\mathscr U_t(c^{-1/2}h)=s^{-1}\mathscr U_t h\), and $s$ is defined in Lemma \ref{l4.1}. In particular,
\[
\mathcal T_{k,t}^{-1}
=\e^{-ik\vartheta(t)}\mathscr U_t^{-1}s.
\]
For a function $u$, a direct calculation yields
\[
\mathscr U_t\partial_y\mathscr U_t^{-1}u
=\rho^{-1/2}\rho\,\partial_Y(\rho^{1/2}u)
=\rho\,\partial_Yu+\tfrac12(\partial_Y\rho) u.
\]
Consequently,
\begin{equation}\label{4.5c}
\mathcal Q_{k,t}
:=\mathscr U_t (-\partial_y^2+|k|^2)\mathscr U_t^{-1}
=-\bigl(\rho\partial_Y+\tfrac12\partial_Y\rho\bigr)^2+|k|^2.
\end{equation}
It follows that the diffusion term transforms according to
\[
\mathcal T_{k,t}(-\partial_y^2+|k|^2) f_k
=s^{-1}\mathcal Q_{k,t}(s g_k)
=(-\partial_Y^2+|k|^2)g_k+(s^{-1}\mathcal Q_{k,t}s-(-\partial_Y^2+|k|^2))g_k.
\]

Since \(V(t,\Psi(t,Y))=\beta(t)+Y^2\), the background transport term becomes
\[
\mathcal T_{k,t}(ikV f_k)=ik\bigl(\beta(t)+Y^2\bigr)g_k.
\]
For the nonlocal term, by \(c(t,\Psi(t,Y))=s^2(t,Y)\) and \eqref{4.5c},
\[
\begin{aligned}
\mathcal T_{k,t}\bigl(c((-\partial_y^2+|k|^2))^{-1} f_k\bigr)
&=s^{-1}s^2\mathcal Q_{k,t}^{-1}(s g_k)=s\mathcal Q_{k,t}^{-1}(s g_k)\\
&=2((-\partial_Y^2+|k|^2))^{-1}g_k+(s\mathcal Q_{k,t}^{-1}s-2(-\partial_Y^2+|k|^2)^{-1})g_k.
\end{aligned}
\]
Moreover, it holds
\[
\mathcal T_{k,t}(\partial_t f_k)
=\partial_tg_k-(\partial_t\mathcal T_{k,t})(\mathcal T_{k,t}^{-1}g_k).
\]
Using
\(\mathcal T_{k,t}=e^{ik\vartheta(t)}\mathscr U_t c^{-1/2}\) yields
\[
(\partial_t\mathcal T_{k,t})\mathcal T_{k,t}^{-1}
=ik\beta(t)
+(\partial_t\mathscr U_t)\mathscr U_t^{-1}
-\frac12\left(\frac{\partial_t c}{c}\right)\circ\Psi.
\]
In the last term, \(\partial_t c(t,y)\) is taken at fixed \(y\) and then evaluated at \(y=\Psi(t,Y)\).

To compute the coordinate contribution, set
\[
b(t,Y):=\frac{\partial_t\Psi(t,Y)}{\partial_Y\Psi(t,Y)}
=\frac{\Psi_t(t,Y)}{J(t,Y)}.
\]
For a fixed function \(h(y)\), write \(u=\mathscr U_t h\). Differentiation at fixed \(Y\) gives
\[
(\partial_t\mathscr U_t)h
=\frac{J_t}{2J}u+J^{1/2}h_y(\Psi)\Psi_t.
\]
The \(Y\)-derivative of \(u=J^{1/2}h(\Psi)\) satisfies
\[
\partial_Yu-\frac{J_Y}{2J}u
=J^{1/2}h_y(\Psi)J.
\]
Substituting this identity into the time derivative and using \(J_t=\partial_Y\Psi_t\), we obtain
\[
\begin{aligned}
(\partial_t\mathscr U_t)h
&=b\,\partial_Yu
+\frac12\left(\frac{J_t}{J}-b\frac{J_Y}{J}\right)u\\
&=b\,\partial_Yu+\frac12(\partial_Yb)u.
\end{aligned}
\]
Hence the time-dependent correction is
\begin{equation}\notag
\mathcal R_t
=b(t,Y)\partial_Y+\frac12\partial_Y b(t,Y)
-\frac12\left(\frac{\partial_t c}{c}\right)\circ\Psi,
\end{equation}
and
\begin{equation}\label{4.6c}
\mathcal T_{k,t}(\partial_t f_k)
=\partial_tg_k-ik\beta(t)g_k-\mathcal R_tg_k.
\end{equation}

Applying \(\mathcal T_{k,t}\) to \eqref{4.1} and using the term \(-ik\beta(t)g_k\) from the time derivative cancels \(+ik\beta(t)g_k\) from the background transport give
\begin{align}\label{4.7c}
&\partial_tg_k
+\nu (-\partial_Y^2+|k|^2)g_k
+ikY^2g_k
+2ik((-\partial_Y^2+|k|^2))^{-1}g_k\notag\\ 
=&\mathcal T_{k,t}\mathrm{div}_k F
-\nu(s^{-1}\mathcal Q_{k,t}s-(-\partial_Y^2+|k|^2))g_k\notag\\ 
&-ik(s\mathcal Q_{k,t}^{-1}s-2(-\partial_Y^2+|k|^2)^{-1})g_k
+\mathcal R_tg_k.
\end{align}

\begin{lemma}\label{l4.2}
Denote \(\|h\|_{H_{k,Y}^{-1}}^2
:=
\|(-\partial_Y^2+|k|^2)^{-\frac{1}{2}}h\|_{L_Y^2}^2.\)
If \(|k|\ge1\) and \(c_*\ll 1\), it holds that
\begin{equation}\label{4.8c}
\|(s^{-1}\mathcal Q_{k,t}s-(-\partial_Y^2+|k|^2))h\|_{H_{k,Y}^{-1}}
\lesssim c_*\|\nabla_kh\|_{L_{Y}^2},
\end{equation}
\begin{equation}\label{4.9c}
\|\nabla_k(s\mathcal Q_{k,t}^{-1}s-2(-\partial_Y^2+|k|^2)^{-1})h\|_{L_{Y}^2}
\lesssim c_*\|h\|_{H_{k,Y}^{-1}},
\end{equation}
\begin{equation}\label{4.9cc}
    \|\mathcal R_t h\|_{H_{k,Y}^{-1}}
\lesssim \nu c_*\|h\|_{L^2_Y}\lesssim \nu c_*\|\nabla_kh\|_{L^2_Y}.
\end{equation}
\end{lemma}
\begin{proof}
For convenience, we denote $Q_k^Y:=(-\partial_Y^2+|k|^2)$. It follows from Lemma \ref{l4.1} that
\[
\|(\rho\partial_Y+\tfrac12\partial_Y\rho-\partial_Y)h\|_{L^2}=\|(\rho-1)\partial_Yh+\tfrac12(\partial_Y\rho) h\|_{L^2}
\lesssim c_*\|\nabla_kh\|_{L^2}.
\]
And the elliptic form
\[
\mathfrak a_{k,t}(h,\varphi)
=
\langle(\rho\partial_Y+\tfrac12\partial_Y\rho)h,(\rho\partial_Y+\tfrac12\partial_Y\rho)\varphi\rangle
+
k^2\langle h,\varphi\rangle
\]
satisfies
\[
\begin{aligned}
\left|\big\langle(\mathcal Q_{k,t}-Q_k^Y)h,\varphi\big\rangle\right|=
\bigl|\mathfrak a_{k,t}(h,\varphi)
-
\langle h_Y,\varphi_Y\rangle
-k^2\langle h,\varphi\rangle
\bigr|
\le
Cc_*
\|\nabla_kh\|_{L^2}
\|\nabla_k\varphi\|_{L^2}.
\end{aligned}
\]
Thus, by Lemma \ref{l4.1}, we have
\begin{equation}\label{4.10c}
  \|(\mathcal Q_{k,t}-Q_k^Y)h\|_{H_{k,Y}^{-1}}
\le
Cc_*\|\nabla_kh\|_{L^2},  
\end{equation}
and for sufficiently small \(c_*\),
\[
\left|\big\langle(\mathcal Q_{k,t}h,h\big\rangle\right|=|\mathfrak a_{k,t}(h,h)|
\ge
\frac12\|\nabla_kh\|_{L^2}^2.
\]
Then if \(\mathcal Q_{k,t}u=r\), testing the equation against \(u\) gives
\[
\frac12\|\nabla_k u\|_{L^2}^2
\le
|\langle r,u\rangle|
\le
\|r\|_{H_{k,Y}^{-1}}\|\nabla_k u\|_{L^2}.
\]
This means that
\[
\|\nabla_k(\mathcal Q_{k,t}^{-1}r)\|_{L^2}
\le
2\|r\|_{H_{k,Y}^{-1}}.
\]
Notice that 
\[
\mathcal Q_{k,t}^{-1}-(Q_k^Y)^{-1}
=
\mathcal Q_{k,t}^{-1}
(Q_k^Y-\mathcal Q_{k,t})
(Q_k^Y)^{-1},
\]
then by \eqref{4.10c}, we get
\begin{equation}\label{4.11c}
\|
\nabla_k((\mathcal Q_{k,t}^{-1}-(Q_k^Y)^{-1})f)
\|_{L^2}
\le Cc_*\|f\|_{H_{k,Y}^{-1}}.
\end{equation}
For the diffusion error, write
\[
s^{-1}\mathcal Q_{k,t}s-(-\partial_Y^2+|k|^2)
=
s^{-1}(\mathcal Q_{k,t}-Q_k^Y)s
+
s^{-1}[Q_k^Y,s].
\]
The first term satisfies the desired estimate by \eqref{4.10c}. For the commutator,
\[
[Q_k^Y,s]h
=
-2s_Yh_Y-s_{YY}h
=
-\partial_Y(s_Yh)-s_Yh_Y.
\]
Then by Lemma \ref{l4.1}, it holds
\[
\|[Q_k^Y,s]h\|_{H_{k,Y}^{-1}}
\le
Cc_*\|\nabla_kh\|_{L^2},
\]
which together with \eqref{4.10c} give
\[
\begin{aligned}
\|(s^{-1}\mathcal Q_{k,t}s-(-\partial_Y^2+|k|^2))h\|_{H_{k,Y}^{-1}}
&\lesssim
\|s^{-1}(\mathcal Q_{\kappa,t}-Q_\kappa^Y)(sh)\|_{H_k^{-1}}
+\|s^{-1}[Q_\kappa^Y,s]h\|_{H_k^{-1}}\\
&\lesssim
c_*\|\nabla_k(sh)\|_{L^2}+c_*\|\nabla_k h\|_{L^2}
\\
&\lesssim c_*\|\nabla_k h\|_{L^2}.
\end{aligned}
\]
This proves \eqref{4.8c}.

For the nonlocal error, put \(s_1=s-\sqrt2\). Then
\[
s\mathcal Q_{k,t}^{-1}s-2(-\partial_Y^2+|k|^2)^{-1}
=
s\bigl(
\mathcal Q_{k,t}^{-1}-(Q_k^Y)^{-1}
\bigr)s+
s_1(Q_k^Y)^{-1}s
+
\sqrt2(Q_k^Y)^{-1}s_1.
\]
Using \eqref{4.11c}, the elliptic estimate for \((Q_k^Y)^{-1}\), and Lemma \ref{l4.1}, we infer
\begin{align}
   &\|\nabla_k(s\mathcal Q_{k,t}^{-1}s-2(-\partial_Y^2+|k|^2)^{-1})h\|_{L^2}\notag\\
   \lesssim&\|\nabla_k(s\bigl(
\mathcal Q_{k,t}^{-1}-(Q_k^Y)^{-1}
\bigr)s)h\|_{L^2}+\|\nabla_k(s_1(Q_k^Y)^{-1}s)h\|_{L^2}+\|\nabla_k((Q_k^Y)^{-1}s_1)h\|_{L^2}\notag\\
\lesssim&\|\nabla_k(\bigl(
\mathcal Q_{k,t}^{-1}-(Q_k^Y)^{-1}
\bigr)s)h\|_{L^2}+c_*\|\nabla_k((Q_k^Y)^{-1}s)h\|_{L^2}+\|s_1h\|_{H_{k,Y}^{-1}}\notag\\
\lesssim&c_*\|h\|_{H_{k,Y}^{-1}}.\notag
\end{align}
This gives \eqref{4.9c}.

Finally, by Lemma \ref{l4.1}, \(\Psi(t,0)=y_c(t)\). For every \(Y\in\mathbb R\), we have 
\[
V(t,\Psi(t,Y))-V(t,\Psi(t,0))=Y^2,
\qquad
V_y(t,\Psi(t,0))=0.
\]
Differentiating the first identity in \(Y\) gives
\[
V_y(t,\Psi(t,Y))J(t,Y)=2Y.
\]
Fixing $Y$ and differentiating with respect to $t$, the time derivative of $Y^2$
  on the right-hand side of the equation is zero, so
  \[
V_t(t,\Psi(t,Y))
   +V_y(t,\Psi(t,Y))\Psi_t(t,Y)
   -V_t(t,\Psi(t,0))=0.
\]
Here we used \(V_y(t,\Psi(t,0))\). Then, combined with the definitions of \(b(t,Y)\) and \(J(t,Y)\), we obtain
\[
\begin{aligned}
2Yb(t,Y)
&=V_y(t,\Psi(t,Y))\partial_t\Psi(t,Y)=V_t(t,\Psi(t,0))-V_t(t,\Psi(t,Y))\\
&=-\nu\bigl[V^L_{yy}(t,\Psi(t,Y))-V^L_{yy}(t,\Psi(t,0))\bigr].
\end{aligned}
\]
Here we used \(V_t=V^L_t=\nu V^L_{yy}\). Therefore, we arrive at
\[
b(t,Y)
=-\frac{\nu}{2}\int_0^1 V^L_{yyy}(t,\Psi(t,sY))J(t,sY)\dd s.
\]
The integral formula extends continuously to \(Y=0\). Differentiating it in \(Y\) yields
\[
b_Y(t,Y)
=-\frac{\nu}{2}\int_0^1 s\Big[
V^L_{yyyy}(t,\Psi(t,sY))J(t,sY)^2
+V^L_{yyy}(t,\Psi(t,sY))J_Y(t,sY)
\Big]\dd s.
\]
The uniform bounds for \(J\) and \(J_Y\) therefore give
\[
\|b\|_{W^{1,\infty}_Y}
\lesssim \nu\bigl(
\|V^L_{yyy}\|_{L^\infty_y}
+\|V^L_{yyyy}\|_{L^\infty_y}
\bigr).
\]
Also, \(c=2+V^L_{yy}\) is bounded below by a positive constant and \(c_t=\nu V^L_{yyyy}\). Hence
\[
\|\frac{c_t}{c}(t,\Psi(t,Y))\|_{L^\infty_Y}
\lesssim \nu\|V^L_{yyyy}\|_{L^\infty_y}.
\]
The heat bound \(\|V^L(t)\|_{H^6}\le\|V_{\mathrm{in}}\|_{H^6}\) and Sobolev embedding now give the asserted coefficient estimate.

Finally, for any \(g\in H_k^1\), we have
\[
\begin{aligned}
|\langle \mathcal R_t h,g\rangle|
&\lesssim
\|b\|_{L^\infty}\|h\|_{L^2}\|\partial_Yg\|_{L^2}
+\bigl(\|b_Y\|_{L^\infty}+\|\frac{c_t}{c}(t,\Psi(t,Y))\|_{L^\infty}\bigr)\|h\|_{L^2}\|g\|_{L^2}\\
&\lesssim \nu\|V_{\mathrm{in}}\|_{H^6}\|h\|_{L^2}\|\nabla_kg\|_{L^2}\lesssim \nu c_*\|h\|_{L^2}\|\nabla_kg\|_{L^2}.
\end{aligned}
\]
This proves \eqref{4.9cc}.
\end{proof}

We now proceed to prove Proposition \ref{p4.1}.
\begin{proof}[Proof of Proposition \ref{p4.1}]
From Lemma \ref{l4.1} and \eqref{4.4c}, one can see that the changes of variables and the multipliers \(c^{\pm1/2}\) satisfy uniform \(H^1\) bounds. In particular,
\begin{equation}\label{4.17a}
\|\nabla_k\mathcal T_{k,t}f\|_{L^2}
\approx
\|\nabla_kf\|_{L^2},
\qquad
\|\mathcal T_{k,t}f\|_{L^2}\approx\|f\|_{L^2}.
\end{equation}

To obtain the negative-order comparison, use duality. Ignoring the scalar phase, which has modulus one,
\[
\mathscr T_t=\mathscr U_tc^{-1/2},
\qquad
\mathscr T_t^*=c^{-1/2}\mathscr U_t^{-1}.
\]
Both \(\mathscr T_t^*\) and \((\mathscr T_t^{-1})^*\) act boundedly on \(H_k^1\). Therefore,
\begin{equation}\label{4.18}
\|\mathcal T_{k,t}f\|_{H_{k,Y}^{-1}}
\approx
\|f\|_{H_{k,y}^{-1}}.
\end{equation}
This comparison identifies the transformed damping norm with the physical velocity norm. If
\[
(\partial_y^2-k^2)\psi_k=f_k,
\qquad
u_{k}=(-\partial_y\psi_k,ik\psi_k),
\]
then
\[
\|u_{k}\|_{L^2}^2
=
\|(-\partial_y^2+|k|^2)^{-1/2}f_k\|_{L^2}^2
=
\|f_k\|_{H_{k,y}^{-1}}^2.
\]
The coordinate transformation therefore preserves the divergence-forcing estimate without losing a horizontal derivative.

Let \(P_\sigma^Y\) be the physical-\(Y\) representation of the bounded multiplier constructed in the preceding section:
\[
P_\sigma^Y
=
\mathcal F_Y^{-1}(I-\delta_0\sigma M)\mathcal F_Y,
\qquad
\sigma=\operatorname{sgn}k,
\]
where \(\delta_0>0\) is the fixed small constant. Set
\[
\mathcal E_\sigma(g)=\langle P_\sigma^Yg,g\rangle.
\]
Similarly to \eqref{3.75}, we have
\[
\mathcal E_\sigma(g_k)\approx\|g_k\|_{L^2}^2,
\]
and
\[
\|\nabla_k P_\sigma^Yg_k\|_{L^2}
\le C\|g_k\|_{H_{k,Y}^1},
\qquad
\|P_\sigma^Yg\|_{H_{k,Y}^{-1}}
\le C\|g\|_{H_{k,Y}^{-1}}.
\]
If $g_k$ solves \eqref{4.7c}, then similarly to \eqref{3.81}, we can prove
\begin{align}\label{4.17}
&\frac{\dd}{\dd t}\mathcal E_\sigma(g_k)
+\nu\|(-\partial_Y^2+|k|^2)^{\frac{1}{2}}g_k\|_{L^2}^2+ c_1|k|\|(-\partial_Y^2+|k|^2)^{-\frac{1}{2}}g_k\|_{L^2}^2+c_1 \nu^{\frac{1}{2}}|k|^{\frac{1}{2}} \|g_k\|_{L^2}^2\notag\\ 
\lesssim&
\operatorname{Re}\langle \mathcal T_{k,t}\mathrm{div}_k F
-\nu(s^{-1}\mathcal Q_{k,t}s-(-\partial_Y^2+|k|^2))g_k\notag\\ 
&\ \ \ \ \ \ \ -ik(s\mathcal Q_{k,t}^{-1}s-2(-\partial_Y^2+|k|^2)^{-1})g_k
+\mathcal R_tg_k,P_\sigma^Yg_k\rangle.
\end{align}
By using Lemma \ref{l4.2}, 
\[
\begin{aligned}
\nu|\langle(s^{-1}\mathcal Q_{k,t}s-(-\partial_Y^2+|k|^2))g,P_\sigma^Yg\rangle|
&\le
\nu\|(s^{-1}\mathcal Q_{k,t}s-(-\partial_Y^2+|k|^2))g\|_{H_{k,Y}^{-1}}
\|\nabla_k P_\sigma^Yg\|_{L^2}\\
&\le
Cc_*\nu\|\nabla_k g\|_{L^2}^2,
\end{aligned}
\]
\begin{align}
&|k||\langle (s\mathcal Q_{k,t}^{-1}s-2(-\partial_Y^2+|k|^2)^{-1})g_k,P_\sigma^Yg_k\rangle|\notag\\
&\le
|k|\|\nabla_k((s\mathcal Q_{k,t}^{-1}s-2(-\partial_Y^2+|k|^2)^{-1})g_k)\|_{L^2}
\|P_\sigma^Yg_k\|_{H_{k,Y}^{-1}}\notag\\
&\le
Cc_*|k|\|g\|_{H_{k,Y}^{-1}}^2,\notag
\end{align}
and
\[
|\langle\mathcal R_tg,P_\sigma^Yg\rangle|
\le
\|\mathcal R_tg\|_{H_k^{-1}}
\|\nabla_k P_\sigma^Yg\|_{L^2}\le
C\nu c_*\|\nabla_kg_k\|_{L^2}^2.
\]
In addition, by using \(0<c_*\ll 1\), 
\[
\|\mathcal T_{k,t}\mathrm{div}_k F\|_{H_{k,Y}^{-1}}
\le C\|F_k\|_{L^2},
\]
and then
\[
\begin{aligned}
2|\langle\mathcal T_{k,t}\mathrm{div}_kF,P_\sigma^Yg_k\rangle|
&\le
C\|F\|_{L^2}\|\nabla_kg_k\|_{L^2}\\
&\le
\frac{1}{4}\nu\|\nabla_kg_k\|_{L^2}^2
+C\nu^{-1}\|F\|_{L^2}^2,
\end{aligned}
\]
we get that
\[
\frac{\dd}{\dd t}\mathcal E_\sigma(g_k)
+c_2\bigl(
\nu\|\nabla_kg_k\|_{L^2}^2
+|k|\|g_k\|_{H_{k,Y}^{-1}}^2
+\nu^{\frac{1}{2}}|k|^{\frac{1}{2}}\|g_k\|_{L^2}^2
\bigr)
\le
C\nu^{-1}\|F_k\|_{L^2}^2.
\]
Therefore, similarly to \eqref{3.73}, by applying \eqref{4.17a}--\eqref{4.18}, and
\[
\mathcal E_\sigma(\mathcal T_{k,t}f_k)\approx\|f_k\|_{L^2}^2,
\]
we can ultimately obtain \eqref{4.2}.
\end{proof}

\section{Nonlinear stability}\label{s5}

In this section, we prove Theorem \ref{thm1}. From \eqref{2.4}, one can see that for $k\in \mathbb{Z} \setminus \{0\}$ and $y\in \mathbb{R}$,
\begin{equation}\label{5.1}
\begin{cases}
\partial_t\omega_k + \nu \Delta_k \omega_k + i k(V + \partial_y^2V\Delta_k^{-1}) \omega_k = -ik(F_1)_k-\partial_y(F_2)_k, \quad\\[2mm] 
U_t - \nu \partial_y^2U = \overline{u_{2,\neq} \omega_{\neq}},\\[2mm]
\partial_t \bar{u}_1- \nu \partial_y^2 \bar{u}_1 = (g_2)_0,\\[2mm]
\omega_k(0,y) = \omega_{k, \mathrm{in}},\ \ \ \ \ U(0,y)=0,
\end{cases}
\end{equation}
where 
\[
(F_1)_k=(g_1)_k-u_{1,k}\partial_yU+\omega_{k}U,\ \ \ (F_2)_k=(g_2)_k-u_{2,k}\partial_yU,
\]
and
\[
(g_1)_k=\sum_{l\neq0,k\neq l} (u_1)_l \omega_{k-l}, \qquad (g_2)_k = \sum_{l\neq0,k\neq l} (u_2)_l \omega_{k-l}.
\]
Recall that \(E_0=\|\overline{\omega}\|_{L^\infty L^2}\) and for \(k \neq 0\),
\begin{align}\label{5.2}
   E_k=& \|\e^{\delta\nu^{\frac{1}{2}}t}\omega_k\|_{L^\infty L^2} + \nu^{\frac{1}{2}}\|\e^{\delta\nu^{\frac{1}{2}}t}(k,\partial_y)\omega_k\|_{L^2 L^2} + (\nu|k|)^{\frac{1}{4}}\|\e^{\delta\nu^{\frac{1}{2}}t}\omega_k\|_{L^2 L^2}  \notag\\
& + |k|^{\frac{1}{2}}\|\e^{\delta\nu^{\frac{1}{2}}t}u_k\|_{L^2 L^2}.
\end{align}

Now we are ready to prove the Theorem \ref{thm1}.

\begin{proof}[Proof of Theorem \ref{thm1}]
The global well-posedness of the 2-D Navier-Stokes equations on the domain $\mathbb{T} \times \mathbb{R} $ around the Poiseuille flow can be obtained by the classical theory for Navier-Stokes equations. Therefore, we only prove the stability estimates.

We first derive the estimate for $E_0$. 
By integration by parts, one has
\[
\left\langle (\partial_t - \nu \partial_y^2)\overline{u}^1, -\partial_y^2 \overline{u}^1 \right\rangle = \frac{1}{2} \partial_t \|\partial_y \overline{u}^1(t)\|_{L^2}^2 + \nu \|\partial_y^2 \overline{u}^1(t)\|_{L^2}^2 = \left\langle (g_2)_0, \partial_y^2 \overline{u}^1 \right\rangle,
\]
which gives
\[
\partial_t \|\partial_y \overline{u}^1\|_{L^2}^2 + \nu \|\partial_y^2 \overline{u}^1\|_{L^2}^2 \lesssim \nu^{-1} \|(g_2)_0\|_{L^2}^2,
\]
this along with \(\partial_y \overline{u}^1(t, y) = \overline{w}(t, y)\) shows
\begin{equation}\label{5.3}
E_0 = \|\overline{\omega}\|_{L^\infty L^2} \lesssim \nu^{-\frac{1}{2}} \|(g_2)_0\|_{L^2 L^2} + \|\overline{\omega}_{\rm in}\|_{L^2}.
\end{equation}

Now we estimate $E_k$. From \eqref{4.2} in Proposition \ref{p4.1} and \eqref{5.2}, we have for $k\neq 0$,
\begin{align}
 E_k\lesssim& \|\omega_{k, \mathrm{in}}\|_{L^2}+ \nu^{-\frac{1}{2}}(\|\e^{\delta \nu^{\frac{1}{2}} t}(g_2)_k\|_{L^2L^2}+\|\e^{\delta \nu^{\frac{1}{2}} t}u_{k}\partial_yU\|_{L^2L^2})\notag\\ 
 &+\nu^{-\frac{1}{2}}(\|\e^{\delta \nu^{\frac{1}{2}} t}(g_1)_k\|_{L^2L^2}+\|\e^{\delta \nu^{\frac{1}{2}} t}\omega_{k}U\|_{L^2L^2}).
\end{align}
For $k\neq 0$, we can derive
\begin{align}\label{5.5}
    \|\e^{\delta \nu^{\frac{1}{2}} t}(g_2)_k\|_{L^2L^2} &\lesssim \sum_{l\neq0,k\neq l} \|\e^{\delta \nu^{\frac{1}{2}} t}(u_2)_l\omega_{k-l}\|_{L^2L^2} \notag\\
&\lesssim \sum_{l\neq0,k\neq l} \|\e^{\delta \nu^{\frac{1}{2}} t}(u_2)_l\|_{L^2L^\infty} \|\e^{\delta \nu^{\frac{1}{2}} t}\omega_{k-l}\|_{L^\infty L^2} \notag\\
&\lesssim \sum_{l\neq0,k\neq l} \|\e^{\delta \nu^{\frac{1}{2}} t}(u_2)_l\|_{L^2L^2}^{\frac{1}{2}} \|\e^{\delta \nu^{\frac{1}{2}} t}\partial_y(u_2)_l\|_{L^2L^2}^{\frac{1}{2}} E_{k-l} \notag\\
&\lesssim \sum_{l\neq0,k\neq l} (|l|^{-\frac{1}{2}}E_l)^{\frac{1}{2}} (|l|\|\e^{\delta \nu^{\frac{1}{2}} t}(u_1)_l\|_{L^2L^2})^{\frac{1}{2}} E_{k-l} \notag\\
&\lesssim \sum_{l\neq0,k\neq l} (|l|^{-\frac{1}{2}}E_l)^{\frac{1}{2}} (|l|^{-\frac{1}{2}}|l|E_l)^{\frac{1}{2}} E_{k-l} \lesssim \sum_{l\neq0,k\neq l} E_lE_{k-l}.
\end{align}
In a similar manner, we can also deduce
\begin{align}\label{5.6}
    \|(g_2)_0\|_{L^2L^2} \lesssim  \sum_{l\neq0} E_lE_{-l}\lesssim(\sum_{l \neq 0}E_l)^2.
\end{align}
Moreover, for $k\neq 0$, there holds
\begin{align}\label{5.7}
  \|\e^{\delta \nu^{\frac{1}{2}} t}(g_1)_k\|_{L^2L^2}=&\sum_{l \in \mathbb{Z} \setminus \{0, k\}}\|\e^{\delta \nu^{\frac{1}{2}} t}(u_1)_l\omega_{k-l}\|_{L^2L^2} \lesssim  \sum_{l \in \mathbb{Z} \setminus \{0, k\}}\| \e^{\delta \nu^{\frac{1}{2}} t} (u_1)_l \|_{L^4 L^\infty} \| \e^{\delta \nu^{\frac{1}{2}} t} \omega_{k-l} \|_{L^4 L^2}\notag\\ 
  \lesssim & \sum_{l \in \mathbb{Z} \setminus \{0, k\}}\| \e^{\delta \nu^{\frac{1}{2}} t} (u_1)_l \|_{L^2 L^2}^{\frac{1}{2}} \|\e^{\delta \nu^{\frac{1}{2}} t} \partial_y(u_1)_l \|_{L^\infty L^2}^{\frac{1}{2}} \| \e^{\delta \nu^{\frac{1}{2}} t} \omega_{k-l} \|_{L^\infty L^2}^{\frac{1}{2}}\| \e^{\delta \nu^{\frac{1}{2}} t} \omega_{k-l} \|_{L^2 L^2}^{\frac{1}{2}}\notag\\
  \lesssim & \nu^{-\frac{1}{8}}\sum_{l \in \mathbb{Z} \setminus \{0, k\}}E_lE_{k-l}.
\end{align}
Next, before estimating the remaining nonlinear terms, we first need to provide the relevant estimates for $U$. From \eqref{5.1}$_3$, it holds
\begin{align}
\frac12\frac{\dd}{\dd t}\|\partial_yU\|_{L^2}^2+\nu\|\partial_y^2U\|_{L^2}^2
=&
{\rm Re}\langle \partial_y(\overline{u_{2,\neq} \omega_{\neq}}),U_y\rangle=-{\rm Re}\langle \overline{u_{2,\neq} \omega_{\neq}},\partial_yU_y\rangle\notag\\ 
\leq&\frac\nu2\|\partial_y^2U\|_{L^2}^2+C\nu^{-1}\|\overline{u_{2,\neq} \omega_{\neq}}\|_{L^2}^2
\end{align}
Since \(\partial_yU(0,y)=0\), by using \eqref{5.6}, we infer
\begin{align}\label{5.9}
    \|\partial_yU\|_{L^\infty L^2}\lesssim \nu^{-\frac{1}{2}}\|\overline{u_{2,\neq} \omega_{\neq}}\|_{L^2L^2}\lesssim \nu^{-\frac{1}{2}}\sum_{ l \neq 0}E_lE_{-l}\lesssim \nu^{-\frac{1}{2}}(\sum_{l \neq 0}E_l)^2.
\end{align}

For \(U\), the heat representation is
\[
U(t)=\int_0^t \e^{\nu(t-s)\partial_y^2}\overline{u_{2,\neq} \omega_{\neq}}(s)\dd s.
\]
The one-dimensional heat-kernel estimate gives
\begin{align}
\|U\|_{L^\infty}
\lesssim&\nu^{-\frac{1}{4}}
\int_0^t(t-s)^{-\frac{1}{4}}\|\overline{u_{2,\neq} \omega_{\neq}}(s)\|_{L^2}\dd s\notag\\
\lesssim&\nu^{-\frac{1}{4}}
\int_0^t(t-s)^{-\frac{1}{4}}\e^{-\delta \nu^{\frac{1}{2}} s}\|\e^{\delta \nu^{\frac{1}{2}} s}\overline{u_{2,\neq} \omega_{\neq}}(s)\|_{L^2}\dd s\notag\\
\lesssim&\nu^{-\frac{1}{4}}\|\e^{\delta \nu^{\frac{1}{2}} t}\overline{u_{2,\neq} \omega_{\neq}}\|_{L^2 L^2}\left(\int_0^t(t-s)^{-\frac{1}{2}}\e^{-2\delta \nu^{\frac{1}{2}} s}\dd s\right)^{\frac{1}{2}}.\notag
\end{align}
For \(t\ge0\), we set \(\tau=(2\delta \nu^{\frac{1}{2}})^{-1}\) and split the integral according to whether \(t-s\le\tau\) or \(t-s>\tau\).

If \(t-s\le\tau\), we use \(\e^{-2\delta \nu^{\frac{1}{2}} s}\le1\) and change variables to \(r=t-s\) to get
\[
\begin{aligned}
\int_{\substack{0\le s\le t\\t-s\le\tau}}
(t-s)^{-\frac{1}{2}}\e^{-2\delta \nu^{\frac{1}{2}} s}\dd s
&\le \int_0^{\min\{t,\tau\}}r^{-\frac{1}{2}}\dd r\\
&\lesssim\sqrt{\min\{t,\tau\}}\\
&\lesssim \sqrt{\tau}\lesssim\nu^{-\frac{1}{4}}.
\end{aligned}
\]
If \(t-s>\tau\), we use \((t-s)^{-\frac{1}{2}}\le\tau^{-\frac{1}{2}}\) to get
\[
\begin{aligned}
\int_{\substack{0\le s\le t\\t-s>\tau}}
(t-s)^{-\frac{1}{2}}\e^{-2\delta \nu^{\frac{1}{2}} s}\dd s
\le \tau^{-\frac{1}{2}}\int_0^\infty \e^{-2\delta \nu^{\frac{1}{2}} s}\dd s\lesssim\frac{\tau^{-\frac{1}{2}}}{\nu^{\frac{1}{2}}}\lesssim\nu^{-\frac{1}{4}}.
\end{aligned}
\]
Combining these two estimates yields
\[
\sup_{t\ge0}\int_0^t(t-s)^{-\frac{1}{2}}\e^{-2\delta \nu^{\frac{1}{2}} s}\dd s \lesssim \nu^{-\frac{1}{4}}.
\]
Thus, using \eqref{5.6} again, we obtain
\begin{align}\label{5.11}
\|U\|_{L^\infty L^\infty}
\lesssim\nu^{-\frac{3}{8}}(\sum_{l \in \mathbb{Z} \setminus \{0\}}E_l)^2.
\end{align}

From \eqref{5.9}--\eqref{5.11} and \eqref{5.2}, we deduce for $k\neq 0$ that
\begin{align}\label{5.11k}
    \|\e^{\delta \nu^{\frac{1}{2}} t}u_{k}\partial_yU\|_{L^2L^2}\lesssim&\|\partial_yU\|_{L^\infty L^2}\|\e^{\delta \nu^{\frac{1}{2}} t}u_{k}\|_{L^2L^\infty}\notag\\ 
    \lesssim&\|\partial_yU\|_{L^\infty L^2}\|\e^{\delta \nu^{\frac{1}{2}} t}u_{k}\|_{L^2L^2}^{\frac{1}{2}}\|\e^{\delta \nu^{\frac{1}{2}} t}\partial_yu_{k}\|_{L^2L^2}^{\frac{1}{2}}\notag\\
    \lesssim&\nu^{-\frac{5}{8}}(\sum_{l \in \mathbb{Z} \setminus \{0\}}E_l)^2E_k,
\end{align}
\begin{align}\label{5.12}
    \|\e^{\delta \nu^{\frac{1}{2}} t}\omega_{k}U\|_{L^2L^2}\lesssim&\|U\|_{L^\infty L^\infty}\|\e^{\delta \nu^{\frac{1}{2}} t}\omega_{k}\|_{L^2L^2}\notag\\ 
    \lesssim&\nu^{-\frac{1}{4}}\|U\|^{\frac{1}{2}}_{L^\infty L^2}\|\partial_yU\|^{\frac{1}{2}}_{L^\infty L^2}E_k\notag\\
    \lesssim&\nu^{-\frac{11}{16}}(\sum_{l \in \mathbb{Z} \setminus \{0\}}E_l)^2E_k.
\end{align}
Combining \eqref{5.3}--\eqref{5.7} and \eqref{5.11k}--\eqref{5.12}, we arrive at

\begin{align}\label{5.13}
    \sum_{k \neq 0  }E_k\lesssim \sum_{k \neq 0  }\|\omega_{{\rm in},k}\|_{L^2}+\nu^{-\frac{5}{8}}\sum_{k \neq 0 }\sum_{l\neq0,k\neq l} E_lE_{k-l}+\nu^{-\frac{19}{16}}(\sum_{l \neq 0}E_l)^2\sum_{k \neq 0 }E_k,
\end{align}
\begin{equation}
    E_0\lesssim \|\overline{\omega}_{\rm in}\|_{L^2}+\nu^{-\frac{1}{2}}(\sum_{k \neq 0}E_k)^2.
\end{equation}
Using \eqref{1.4}, one gets that
\[
\sum_{k\in\mathbb Z}\|\omega_{\mathrm{in},k}\|_{L_y^2}
\le
\left(\sum_{k\in\mathbb Z}\frac1{1+k^2}\right)^{1/2}
\left(\sum_{k\in\mathbb Z}(1+k^2)
               \|\omega_{\mathrm{in},k}\|_{L_y^2}^2\right)^{1/2}
\leq C\varepsilon_*\nu^{-\frac{5}{8}}.
\]
Then by the smallness of $\varepsilon_*$, we get \eqref{1.5g}, and hence finish the proof of Theorem \ref{thm1}.

\end{proof}

\noindent{\bf Conflict of Interest:} The authors declare that they have no conflict of interest.

\vspace{0.5cm}
\noindent{\bf
Data Availability Statement} No data were used for the research described in this article.

\vspace{0.5cm}
\noindent{\bf Declaration of Generative AI and AI-Assisted Technologies:} During the preparation of this work, the author used ChatGPT in order to improve language and readability. After using this tool/service, the author reviewed and edited the content as needed and takes full responsibility for the content of the publication.

\vspace{0.5cm}
\noindent{\bf Acknowledgments:} The author welcomes any constructive comments on this manuscript. This research received no funding.

\end{document}